\documentclass[12pt,leqno]{article}
\usepackage{amsfonts}
\usepackage{amsmath, amsthm, amsfonts, amssymb, color}
\usepackage{mathrsfs}
\usepackage{url}
\usepackage{color}
\theoremstyle{definition}

\newcommand{\scr}[1]{\mathscr #1}
\definecolor{wco}{rgb}{0.5,0.2,0.3}

\numberwithin{equation}{section} \theoremstyle{remark}

\newcommand{\ua}{\uparrow}

\title{{\bf Exponential Ergodicity for Singular Reflecting  McKean-Vlasov SDEs  }\footnote{Supported in
 part by 
 NNSFC (11771326, 11831014, 11921001) and the DFG through CRC 1283.} }
\author{{\bf    Feng-Yu Wang  }\\
\footnotesize{Center for Applied Mathematics, Tianjin University, Tianjin 300072, China}\\
\footnotesize{wangfy@tju.edu.cn } }

\begin{document}
\allowdisplaybreaks
\def\R{\mathbb R}  \def\ff{\frac} \def\ss{\sqrt} \def\B{\mathbf
B}\def\TO{\mathbb T}
\def\I{\mathbb I_{\pp M}}\def\p<{\preceq}
\def\N{\mathbb N} \def\kk{\kappa} \def\m{{\bf m}}
\def\ee{\varepsilon}\def\ddd{D^*}
\def\dd{\delta} \def\DD{\Delta} \def\vv{\varepsilon} \def\rr{\rho}
\def\<{\langle} \def\>{\rangle} \def\GG{\Gamma} \def\gg{\gamma}
  \def\nn{\nabla} \def\pp{\partial} \def\E{\mathbb E}
\def\d{\text{\rm{d}}} \def\bb{\beta} \def\aa{\alpha} \def\D{\scr D}
  \def\si{\sigma} \def\ess{\text{\rm{ess}}}
\def\beg{\begin} \def\beq{\begin{equation}}  \def\F{\scr F}
\def\Ric{{\rm Ric}} \def\Hess{\text{\rm{Hess}}}
\def\e{\text{\rm{e}}} \def\ua{\underline a} \def\OO{\Omega}  \def\oo{\omega}
 \def\tt{\tilde}
\def\cut{\text{\rm{cut}}} \def\P{\mathbb P} \def\ifn{I_n(f^{\bigotimes n})}
\def\C{\scr C}      \def\aaa{\mathbf{r}}     \def\r{r}
\def\gap{\text{\rm{gap}}} \def\prr{\pi_{{\bf m},\varrho}}  \def\r{\mathbf r}
\def\Z{\mathbb Z} \def\vrr{\varrho} \def\ll{\lambda}
\def\L{\scr L}\def\Tt{\tt} \def\TT{\tt}\def\II{\mathbb I}
\def\i{{\rm in}}\def\Sect{{\rm Sect}}  \def\H{\mathbb H}
\def\M{\scr M}\def\Q{\mathbb Q} \def\texto{\text{o}} \def\LL{\Lambda}
\def\Rank{{\rm Rank}} \def\B{\scr B} \def\i{{\rm i}} \def\HR{\hat{\R}^d}
\def\to{\rightarrow}\def\l{\ell}\def\iint{\int}
\def\EE{\scr E}\def\Cut{{\rm Cut}}\def\W{\mathbb W}
\def\A{\scr A} \def\Lip{{\rm Lip}}\def\S{\scr S}
\def\BB{\scr B}\def\Ent{{\rm Ent}} \def\i{{\rm i}}\def\itparallel{{\it\parallel}}
\def\g{{\mathbf g}}\def\Sect{{\mathcal Sec}}\def\T{\mathcal T}\def\V{{\mathbb V}}
\def\PP{{\bf P}}\def\HL{{\bf L}}\def\Id{{\rm Id}}\def\f{{\bf f}}\def\cut{{\rm cut}}
\def\n{{\mathbf n}}
\def\sm{\preceq}
\def\BL{\scr A} \def\I{{\mathbf I}}

\maketitle

\begin{abstract} By refining a recent result of  Xie and Zhang \cite{XZ}, we prove the exponential ergodicity under a weighted variation norm for singular
SDEs with drift containing a local integrable term and a coercive term.
This result is then extended to singular reflecting SDEs as well as singular McKean-Vlasov SDEs with or without reflection.
The exponential  ergodicity  in   the relative entropy and (weighted) Wasserstein distances are also studied  for reflecting McKean-Vlasov SDEs.    The main results are illustrated by non-symmetric singular granular media equations.
\end{abstract} \noindent
 \vskip 2cm

\section{Introduction}

Let $D\subset \R^d$ be a connected open domain including the global situation $D=\R^d$, and let $\scr P$ denote the space of probability measures on $\bar D$, the closure of $D$.
Consider the following distribution dependent (i.e. McKean-Vlasov) SDE   on $\bar D$ with reflection if $D\ne \R^d$:
\beq\label{E1} \d X_t= b(X_t,\L_{X_t}) \d t+ \si(X_t)\d W_t + \n(X_t) \d l_t,\ \ t\ge 0,\end{equation}
 where    $(W_t)_{t\ge 0}$ is an $m$-dimensional Brownian motion on a complete filtration probability space $(\OO,\{\F_t\}_{t\ge 0},\P)$,
 $\L_{X_t}$ is the distribution of $X_t$,
$$b: D\times \scr P \to \R^d,\ \ \   \si: D \to \R^d\otimes\R^m$$ are measurable,
and    when $D\ne \R^d$,
 $\n$ is the inward unit normal vector field of the boundary $\pp D$,  and
   $l_t$ is an adapted continuous increasing process which increases only when $X_t\in \pp D$.

In the case that   $D=\R^d$, we have   $l_t=0$ so that  \eqref{E1} becomes
the distribution dependent SDE (DDSDE)
\beq\label{E1b} \d X_t= b(X_t,\L_{X_t}) \d t+ \si(X_t)\d W_t,\ \ t\ge 0.\end{equation}
If moreover $b(x,\mu)=b(x)$ does not depend on $\mu$, it reduces to the classical  It\^o's SDE
\beq\label{E1'} \d X_t = b(X_t)\d t+\si(X_t)\d W_t,\ \ t\ge 0.\end{equation}

In the recent work \cite{W21b}, the well-posedness and regularity estimates have been studied for solutions to \eqref{E1} with $b$ containing a locally integrable term and a Lipchitz continuous term.
However, the ergodicity was only investigated under monotone or Lyapunov conditions  excluding this singular situation.
See also \cite{BUT, EGZ, GLW, HSS, LMW, 20RWb, V06, W21a} and references within for results on the ergodicity  of McKean-Vlasov SDEs without reflection under monotone or Lyapunov conditions.
On the other hand, by using Zvokin's transform, the exponential ergodicity was proved by Xie and Zhang  \cite{XZ} for the singular SDE \eqref{E1'}.
In this paper, we aim to refine the result of   \cite{XZ} and make extensions to singular SDEs with reflection  and   distribution dependent drift.

When the SDE \eqref{E1}  is well-posed, let  $P_t^*\nu=\L_{X_t}$ for the solution with  initial distribution $\nu\in \scr P.$
We will study the exponential convergence  of $P_t^*$  under the weighted variation distance induced by   a positive measurable function $V$:
$$\|\mu-\nu\|_{V} :=|\mu-\nu|(V)=\sup_{|f|\le V}|\mu(f)-\nu(f)|,\ \ \mu,\nu\in \scr P,$$
where $|\mu-\nu|$ is the total variation of $\mu-\nu$ and $\mu(f):=\int f\d\mu$ for a measure $\mu$ and $f\in L^1(\mu)$.
When $V=1$,  $\|\cdot\|_V$ reduces to the the total  variation norm $\|\cdot\|_{var}.$

We will consider   $b(x,\mu)= b^{(0)}(x)+b^{(1)}(x,\mu)$, where $b^{(0)}$ is the singular term satisfying
\beq\label{0*1} \sup_{z\in\R^d} \int_{B(z,1)\cap D} |b^{(0)}(x)|^p(\d x) <\infty\end{equation} for some $p>d\lor 2$,
and  $b^{(1)}(\cdot,\mu) $ is a coercive term  such that
$$\limsup_{x\in\bar D,|x|\to\infty} \sup_{\mu\in \scr P}\<b^{(1)}(x,\mu),\nn V(x)\>=-\infty$$
holds for some compact function $V\in C^2(\R^d)$ (i.e.   $\{V\le r\}$ is compact for any $r>0$). The later condition is trivial for bounded $D$ by taking $V=1$ and the convention that
$\sup\emptyset=-\infty$.

\

 To conclude this section, we present below  an example for the $L^1$-exponential convergence of non-symmetric singular granular media equations,
 see \cite{CMV,GLW,20RWb} for the study  of regular and symmetric models for $D=\R^d$.

\paragraph{Example 1.1.} Let $D=\R^d$ or be a  $C^{2,L}$-domain (see Definition \ref{D2} below).  Consider the following nonlinear PDE for probability density functions on $\bar D$:
\beq\label{GRN}   \pp_t  \varrho_t= \DD  \varrho_t - {\rm div} \big\{ \varrho_t b  +  \varrho_t (W* \varrho_t)\big\}, \ \ \nn_\n \varrho_t|_{\pp D}=0 \ \text{if}\ \pp D\ne\emptyset,\end{equation}
 where \beg{enumerate}
 \item[(i)] $W$ is a bounded measurable function on $\bar D\times \bar D$, and
  $$(W* \varrho_t)(x):=\int_{\R^d} W(x, z)  \varrho_t(z)\d z;$$
\item[(ii)]    $b=b^{(0)}+b^{(1)}$ is a vector field such that  \eqref{0*1} holds for some $p>d\lor 2$, and
$b^{(1)}$ is locally bounded with $b^{(1)}(x)=- \phi(|x|^2) x$ for larger $|x|$ and some increasing function $\phi: [0,\infty)\to [1,\infty)$ with $\int_1^\infty\ff{\d s}{s\phi(s)}<\infty$.\end{enumerate}

 In physics,  $\rr_t$ stands for the distribution density of particles, $W$  describes the interaction among particles,  and $b$ refers to the potential of individual  particles.
 When $b$ and $W$ are not of gradient type, the associated mean field particle systems are non-symmetric.

To characterize \eqref{GRN} using \eqref{E1}, let
 $$b(x,\mu)= b(x) + (W*\mu)(x) ,\ \ \si(x)= \ss 2 {\bf I}_d,$$ where ${\bf I}_d$ is the $d\times d$ identity matrix, and $(W*\mu)(x):= \int_{\bar D} W(x,z)\mu(\d z).$

 By (i) and (ii), {\bf (A1)} holds for $V(x):= |x|^2$ when $D=\R^d$, while {\bf (A2)} holds for $V=1$ when $D$ is a bounded $C_b^{2,L}$ domain. So, by Theorem \ref{T3},  \eqref{E1} is well-posed, and
by It\^o's formula,    $\rr_t(x) := \ff{\d P_t^*\nu}{\d x} $ solves \eqref{GRN} for $\rr_0(x):= \ff{\d\nu}{\d x}$, see Subsection 1.2 in \cite{W21b}.
On the other hand,  when $D=\R^d$ the superposition principle in \cite{RWBR} says that a solution of \eqref{GRN} is the distribution density of a weak solution to \eqref{E1}, such that \eqref{GRN} is well-posed as well. Moreover:
\beg{enumerate}
\item[(a)] By    Theorem \ref{T3}, when $\|W\|_\infty$ is small enough,  $P_t^*$ has a  unique invariant probability measure $\mu$ satisfying \eqref{Q2}, so that the solution $\rr_t:=\ff{\d P_t^*\nu}{\d x}$ of \eqref{GRN}
satisfies
$$\|\rr_t-\rr\|_{L^1} =\|P_t^*\nu-\mu\|_{var} \le c\e^{-\ll t}\|\rr_0-\rr\|_{L^1},\ \ t\ge 0$$
for some constants $c,\ll>0$, where $\rr$ is the density function of $\mu$.
\item[(b)] Let $D=\R^d$ or $D$ be convex.  If there exists a constant $K>0$ such that
\beq\label{*0} \<b(x)-b(y),x-y\>\le -K|x-y|^2,\ \ x,y\in D\end{equation} holds, by Theorem \ref{T6}, when $\|\nn^2 W\|_\infty$ is small enough $P_t^*$ is exponential ergodic
in the relative entropy and the quadratic Wasserstein distance $\W_2$.  If \eqref{*0} only holds for large $|x-y|$,
according to Theorem \ref{T7},  $P_t^*$ is exponential ergodic under a weighted Wasserstein distance provided $\|\nn^2 W\|_\infty$ is small enough.
 \end{enumerate}

\

In the remainder of the paper,   we study in Section 2 the exponential ergodicity for singular reflecting SDEs, then prove the uniform ergodicity for singular reflecting McKean-Vlasov SDEs in Section 3, and finally  investigate  in Section 4 the exponential ergodicity  for reflecting McKean-Vlasov SDEs in relative entropy and (weighted) Wasserstein distances.

\section{Exponential ergodicity for singular reflecting SDEs  }

To measure the singularity of the SDE, we introduce some functional spaces used in \cite{XXZZ}. For any $p\ge 1$, let $L^p$ be the class of   measurable  functions $f$ on $  D$ such that
 $$\|f\|_{L^p}:=\bigg(\int_{D }|f(x)|^p\d x\bigg)^{\ff 1 p}<\infty.$$
For any $\epsilon >0$ and $p\ge 1$, let  $H^{\epsilon,p}:=(1-\DD)^{-\ff\epsilon 2} L^p$ with
 $$\|f\|_{H^{\epsilon,p}}:= \|(1-\DD)^{\ff\epsilon 2} f\|_{L^p}<\infty,\ \ f\in  H^{\epsilon,p},$$  where $\DD$ is the (Neumann if $\pp D\ne\emptyset$) Laplacian.
  For any $z\in\R^d$ and $r>0$,   let
  $$B(z,r):=\{x\in\R^d: |x-z|\le  r\}$$ be the closed ball centered at $z$ with radius $r$. We will simply denote $B_r=B(0,r)$ for $r>0$.
We write $f\in \tt L^p$ if
$$\|f\|_{\tt L^p}:=\sup_{z\in\bar D} \|1_{B(z,1)} f\|_{L^p}<\infty.$$
Moreover, let $g\in C_0^\infty(\bar D)$ with $g|_{B_1}=1$ and the Neumann boundary condition $\nn_\n g|_{\pp D}=0$ if $\pp D$ exists. We denote $f\in \tt H^{\epsilon,p}$ if
$$\|f\|_{\tt H^{\epsilon, p}}:= \sup_{z\in \bar D} \|g(z+\cdot) f\|_{H^{\epsilon,p}}<\infty.$$
We note that  the space $\tt H^{\epsilon,p}$ does not depend on the choice of $g$. If  a vector or matrix valued function has components   in one of the above introduced spaces,
then   it is said in the same space with  norm defined as the sum of  components' norms.

In the following subsections, we first  state the main results, then present some lemmas, and finally prove the main results.

\subsection{Main results}

We  first consider  the ergodicity of  SDE \eqref{E1'} under   the following assumption, where by the Sobolev embedding theorem $\si$ (hence $\si\si^*$) is H\"older continuous by the boundedness of $\si$ and $\|\nn \si\|\in \tt L^p$ for some $p>d$.

 \beg{enumerate} \item[{\bf (A1)}] $\si$ is weakly differentiable, $\si\si^*$ is invertible, and   $b=b^{(0)}+b^{(1)}$ such that  the following conditions hold.
 \enumerate
 \item[$(1)$]  There exists $p>d\lor 2$ such that
 $$\|\si\|_\infty+\|(\si\si^*)^{-1}\|_\infty +\|b^{(0)}\|_{\tt L^p}+ \|\nn \si\|_{\tt L^p}<\infty.$$
 \item[$(2)$]   $b^{(1)}$ is locally bounded, there exist   constants $K>0,\vv\in (0,1)$, some compact function  $V\in C^2(\R^d; [1,\infty))$, and a continuous  increasing function
 $\Phi: [1,\infty)\to [1,\infty)$ with $\Phi(n)\to \infty$ as $n\to\infty$,  such that
 \beq\label{LYP}\beg{split}&\<b^{(1)},\nn V\>(x)+\vv |b^{(1)}(x)| \sup_{B(x,\vv)}  \{|\nn V|+ |\nn^2 V|\} \le K-\vv (\Phi\circ V)(x),\\
 &\lim_{|x|\to\infty}  \ff{\sup_{B(x,\vv)}\{\|\nn^2 V\|+|\nn V|\}} {V(x)\land (\Phi\circ V)(x)}=0.\end{split}\end{equation}
 \end{enumerate}

\beg{thm}\label{T1.6.1} Assume {\bf (A1)}. Then $\eqref{E1'}$ is well-posed,   the associated Markov semigroup $P_t$ has a unique invariant probability measure $\mu$ such that $\mu(\Phi(\vv_0V))<\infty$ for some $\vv_0\in (0,1)$, and
\beq\label{EGD} \lim_{t\to\infty} \|P_t^*\nu-\mu\|_{var} =0,\ \ \nu\in \scr P. \end{equation} Moreover:
\beg{enumerate} \item[$(1)$] If $\Phi(r)\ge \dd r$ for some constant $\dd>0$ and all $r\ge 0$, then there exist constants $c>1,\ll>0$
such that
 \beq\label{EX1}  \|P_t^*\mu_1-P_t^*\mu_2\|_{V}  \le c\e^{-\ll t} \|\mu_1-\mu_2\|_V,\ \ \mu_1,\mu_2\in \scr P, t\ge 0.\end{equation}
 In particular,
$$  \|P_t^*\nu-\mu \|_{V}  \le c\e^{-\ll t} \|\nu-\mu\|_V,\ \ \nu \in \scr P, t\ge 0. $$
 \item[$(2)$] Let $ H(r):=\int_0^r   \ff {\d s} {\Phi(s)} <\infty$ for $r\ge 0$.  If $\Phi$ is convex, then there exist constants $k>1,\ll>0$ such that
 \beq\label{EX0}  \|P_{t}^*\dd_x-\mu\|_{V}\le k \big\{1+H^{-1} (H(V(x))- k^{-1} t)\big\} \e^{-\ll t},\ \ x\in \R^d, t\ge 0,\end{equation}
 where $H^{-1}$ is the inverse of $H$ with $H^{-1}(r):=0$ for $r\le 0$.
  Consequently, if $H(\infty)<\infty$ then there exist constants $c,\ll,t^*>0$ such that
  \beq\label{EX2}\|P_t^*\mu_1-\mu_2\|_{V}\le c  \e^{-\ll t}\|\mu_1-\mu_2\|_{var},\ \ t\ge t^*, \mu_1,\mu_2\in   \scr P.\end{equation} \end{enumerate}
\end{thm}

To illustrate this result, we present below a consequence which covers the  situation of \cite[Theorem 2.10]{XZ} where
$$\<b^{(1)}(x),x\>\le c_1 -c_2 |x|^{1+p},\ \ |b^{(1)}(x)|\le c_1(1+|x|)^{p}$$ holds for some constants $c_1,c_2>0$ and $p\ge 1$.
Indeed, Corollary \ref{C1.2} implies the exponential ergodicity under the weaker condition
\beq\label{B!!} \<b^{(1)}(x),x\>\le c_1 -c_2 |x|^{1+p},\ \ |b^{(1)}(x)|\le c_1(1+|x|)^{p+1}\end{equation}   for some constants $p, c_1,c_2>0$ ($p$ may be smaller than $1$,  $|b^{(1)}|$ may have higher order growth),
since in this case, \eqref{B1} and \eqref{PII} hold for $\phi(r):= (1+ r)^{\ff{1+p}2}$,  and \eqref{PSI} holds for $\psi(r):= (1+r^2)^q$ for any $q>0$ when $p\ge 1.$

 \beg{cor}\label{C1.2} Assume {\bf (A1)}$(1)$ and let $b^{(1)}$ satisfy
 \beq\label{B1} \<b^{(1)}(x),x\>\le c_1-c_2 \phi(|x|^2),\ \ |b^{(1)}(x)|\le c_1 \phi(|x|^2),\ \ x\in \R^d\end{equation}
 for some constants $c_1,c_2>0$ and increasing function $\phi: [0,\infty)\to [1,\infty)$ with
 \beq\label{PII} \aa:= \liminf_{r\to\infty}\ff{\log \phi(r)} {\log r}>\ff 1 2.\end{equation}   Then
  \beg{enumerate}
 \item[$(1)$]   $\eqref{E1'}$ is well-posed,  $P_t$ has a unique invariant probability measure $\mu$ such that
   $\mu(V)<\infty$ and $\eqref{EX1}$  hold for      $V:=\e^{(1+|\cdot|^2)^\theta}$ with $\theta\in ((1-\aa)^+,\ff 1 2)$.
   In general,   for any increasing function $1\le \psi\in C^2([1,\infty))$ satisfying
 \beq\label{PSI} \liminf_{r\to\infty} \ff{\psi'(r)\phi(r)}{\psi(r)}>0,\ \ \lim_{r\to\infty} \ff{\psi''(r)r}{\psi(r)}=0,\end{equation}   $\mu(V)<\infty$ and  $\eqref{EX1}$
hold  for  $V :=\psi(|\cdot|^2)$.
 \item[$(2)$] If $\int_0^\infty\ff{\d s}{\phi(s)}<\infty$, then  $\eqref{EX2}$ holds   $V:=(1+|\cdot|^2)^q (q>0)$ and some constants $c,\ll,t^*>0.$ \end{enumerate}
 \end{cor}

 \paragraph{Remark 2.1.}  We have the following assertions on the invariant probability measure $\mu$ and the ergodicity in Wasserstein distance and relative entropy.
 \beg{enumerate}
  \item[(1)] According to \cite[Corollary 1.6.7 and Theorem 3.4.2]{BKRS},   {\bf (A1)} implies that   $\mu$ has a strictly positive density
 function $\rr\in H_{loc}^{1,p}$, the space of functions $f$ such that $fg\in H^{1,2}$ for all $g\in C_0^\infty(\R^d).$   Moreover, by \cite[Theorem 3.1.2]{BKRS},
 when $\si$ is Lipschitz continuous and $\mu(|b|^2)<\infty,$ we have $\ss\rr\in H^{1,2}.$ So, when $\eqref{B1}$ holds for $\phi(r)\sim r^p$ for some $p>\ff 1 2$ and large $r>0$,
 Corollary \ref{C1.2}(1) implies that $\mu$ has density with $\ss \rr\in H^{1,2}.$
  See also \cite{W17} and \cite{W18}  for  different type global regularity estimates on $\rr$ under
 integrability conditions.
 \item[(2)] Let $V:= (1+|\cdot|^2)^{\ff p 2}$ for some $p\ge 1$. By \cite[Theorem 6.15]{VV},   there exists a constant $c(p)>0$ such that
$$\W_p(\mu,\nu)^p \le c(p) \|\mu-\nu\|_{V},$$  where
$$ \W_p(\mu_1,\mu_2):= \inf_{\pi\in \C(\mu_1,\mu_2)}\bigg( \int_{\R^d\times\R^d} |x-y|^p\pi(\d x,\d y) \bigg)^{\ff 1 p}   $$
for  $\C(\mu_1,\mu_2)$ being  the set of couplings for $\mu_1$ and $\mu_2$.
So,  by Corollary \ref{C1.2},  if {\bf (A1)} holds with   $\Phi(r)\ge \dd r$ for some   $\dd>0, $
then there exist  constants $c,\ll>0$ such that
$$\W_p(P_t^*\nu,\mu)^p\le c(1+\nu(|\cdot|^p) ) \e^{-\ll t},\ \ t\ge 0, \nu\in\scr P;$$  and if  moreover  $\Phi$ is convex with $\int_0^\infty\ff{\d s}{\Phi(s)}<\infty$,
 then  there exist constants $c,\ll,t^*>0$ such that
 $$\W_p(P_t^*\nu,\mu)^p\le c   \e^{-\ll t}\|\mu-\nu\|_{var},\ \ t\ge t^*, \nu\in\scr P.$$
\item[(3)]  When   $b^{(1)}$ is Lipschitz continuous, the log-Harnack inequality in  \cite[Theorem 4.1]{YZ0} implies
$$\Ent(P_t^*\nu|\mu)\le \ff{c'}{1\land t} \W_2(\nu,\mu)^2,\ \ \nu\in \scr P,t>0$$ for some constant $c'>0$,
where $\Ent(\nu|\mu)$ is the relative entropy. Thus, by Corollary \ref{C1.2},  if {\bf (A1)} holds for $V(x):= 1+ |x|^2$ and   $\Phi(r)\ge \dd r$ for some constant $\dd>0,$
then  there exist  constants $c,\ll>0$ such that
$$\Ent(P_t^*\nu|\mu)\le c(1+\nu(|\cdot|^2) ) \e^{-\ll t},\ \ t\ge  1, \nu\in\scr P;$$
 and  if   moreover   $\Phi$ is convex with $\int_0^\infty\ff{\d s}{\Phi(s)}<\infty$, then there exist $c,\ll,t^*>0$ such that
 $$\Ent(P_t^*\nu|\mu)\le c \e^{-\ll t}\|\mu-\nu\|_{var},\ \ t\ge  t^*, \nu\in\scr P.$$ \end{enumerate}

 Next,
consider  the following reflecting SDE on $D\ne \R^d$:
\beq\label{E1'b} \d X_t = b(X_t)\d t+\si(X_t)\d W_t+\n(X_t)\d l_t,\ \ t\ge 0,\end{equation}
where $\pp D\in C_b^{2,L}$ which is defined as follows.

\beg{defn}\label{D2} Let $\rr_\pp$ be the distance function to $\pp D$.
For any  $k\in\mathbb N$, we write $\pp D\in C^k_{b}$  if there exists a constant $r_0>0$ such that  the polar coordinate around  $\pp D$
$$ \pp D\times [-r_0,r_0] \ni (\theta, r)\mapsto   \theta+ r\n(\theta) \in B_{r_0}(\pp D):=\{x\in \R^d: \rr_\pp(x)\le r_0\}$$
is a $C^k$-diffeomorphism. We write $\pp D\in C_b^{k,L}$, if it is $C_b^k$ with $\nn^k\rr_\pp$ being Lipschitz continuous on $B_{r_0}(\pp D).$ \end{defn}

We also need heat kernel estimates for the Neumann semigroup
$\{P_{t}^{\si}\}_{t\ge 0}$  generated by
$$   L^{\si} := \ff 1 2  {\rm tr}\big( \si_t\si_t^* \nn^2\big).$$
For any $\varphi\in C_b^2(\bar D)$,   let $ P_{t}^{\si}\varphi $ be the solution of  the PDE
\beq\label{NMM} \pp_t u_t=  L^{\si} u_t,\ \ \nn_{\n}u_t|_{\pp D}=0\ \text{for}\ s>0, u_0=\varphi.\end{equation}
We will prove the exponential ergodicity of \eqref{E1'b}  under the following assumption.

 \beg{enumerate}  \item[{\bf (A2)}]    $\pp D\in C_b^{2,L}$ and the following conditions hold.
\item[$(1)$]  {\bf (A1)} holds for $\bar D$ replacing $\R^d$, and there exists $r_0>0$ such that
\beq\label{BDR}  \nn_{\n(x)} V(y)\le 0, \ \ x\in \pp D, |y-x|\le r_0. \end{equation}
\item[$(2)$]   For any  $\varphi\in C_b^2(\bar D)$, the PDE \eqref{NMM} has a unique solution $P_{t}^{\si}\varphi\in C_b^{1,2}(\bar D),$  such that    for some  constant $c>0$  we have
$$ \|\nn^i P_{t}^{\si} \varphi\|_{\infty} \le  c (1\land t)^{-\ff 1 2} \|\nn^{i-1} \varphi\|_{\infty},\ \ t>0, \ i=1,2,\varphi\in C_b^2(\bar D),$$  where $\nn^{0}\varphi:=\varphi.$
 \end{enumerate}

As explained in \cite[Remark 2.2(2)]{W21b} that, {\bf (A2)}(2) holds if    $D$ is  bounded  and $\si$ is H\"older continuous.
  Moreover, \eqref{BDR} is trivial when   $\pp D$ is bounded, since in this case we may take $1\le \tt V\in  C^2(\R^d)$ such that $\tt V=1$ on $\pp_{r_0}(\pp D)$ and $\tt V=V$ outside a compact set,
  so that \eqref{LYP} remains true for $\tt V$ replacing $V$.  Similarly,   \eqref{BDR}  holds for $V(x_1,x_2):= V_1(x_1)+V_2(x_2)$ and $D= D_1\times \R^l$ where $l\in \mathbb N$ is less than $d$,
   $\pp D_1\subset \R^{d-l}$ is bounded,  and $V_1=1$ in a neighborhood of $\pp D_1$.

 \beg{thm}\label{T2} Assume {\bf (A2)}. Then all assertions in Theorem \ref{T1.6.1} hold for the reflecting SDE $\eqref{E1'b}$.
  \end{thm}

  \subsection{Some lemmas}

  We first consider the following time dependent SDE with reflection when $\pp D$ exists:
\beq\label{E1.1}   \d X_t = b_t(X_t)\d t+\si_t(X_t)\d W_t+\n(X_t)\d l_t,\ \ t\ge 0.\end{equation}

For any $T>0$ and $p,q>1$, let $\tt L_q^p(T)$ denote the class of measurable functions $f$ on $[0,T]\times\bar D$ such that
$$\|f\|_{\tt L_q^p(T)}:= \sup_{z\in \bar D} \bigg( \int_{0}^{T} \|1_{B(z,1)}f_t\|_{L^p}^q\d t\bigg)^{\ff 1 q}<\infty.$$
For any $\epsilon>0$, let $\tt H_{q}^{\epsilon,p}(T)$ be the space of $f\in \tt L_q^p$ with
$$\|f\|_{\tt H_q^{\epsilon,p}(T)}:=\sup_{z\in \bar D}  \bigg( \int_{0}^{T}  \|f_t\|_{\H^{\epsilon,p}}^q\d t\bigg)^{\ff 1 q}<\infty.$$

We will study the well-posedness, strong Feller property and irreducibility under the following assumptions for $D=\R^d$ and $D\ne \R^d$ respectively.

\beg{enumerate}
\item[{\bf (A3)}]        Let   $T>0, D=\R^d$,  $a_t(x):= (\si_t\si^*_t)(x)$ and $b_t(x)= b_t^{(0)}(x)+ b_t^{(1)}(x)$.
\item[$(1)$] $a$  is invertible with $\|a\|_\infty+\|a^{-1}\|_\infty<\infty$    and
$$ \lim_{\vv\to 0} \sup_{|x-y|\le \vv, t\in [0,T]} \|a_t(x)-a_t(y)\|=0.$$
\item[$(2)$] There exist $l\ge 1$,  $\{ (p_i,q_i)\}_{0\le i\le l}\in \scr K:=\{(p,q): p,q\in (2,\infty), \ff d p+\ff 2 q<1\}$  and $1\le f_i\in \tt L_{q_i}^{p_i}$   such that
$$  |b^{(0)}|\le f_0,\ \ \|\nn\si\|\le \sum_{i=1}^l f_i.$$
\item[$(3)$] $b^{(1)}$ is locally bounded,   there exist constants $K,\vv>0$, increasing $\phi\in C^1([0,\infty); [1,\infty))$ with $\int_0^\infty\ff{\d s}{r+\phi(s)}=\infty$,
and  a compact function  $V\in C^2(\R^d; [1,\infty))$ such that
\beg{align*} &\sup_{B(x, \vv)}\big\{|\nn V|+   \|\nn^2 V\| \big\}\le K V(x),\\
&   \<b_t^{(1)}(x),\nn V(x) \>+\vv |b_t^{(1)}(x)|\sup_{B(x,\vv)} \|\nn^2 V\| \le K\phi(V(x)),  \ \ (t,x)\in  [0,T]\times \R^d.\end{align*} \end{enumerate}

When  $D\ne \R^d$, we consider the following time dependent differential operator  on $\bar D$:
\beq\label{BARL}    L_t^{\si} := \ff 1 2  {\rm tr}\big( \si_t\si_t^* \nn^2\big),\ \ t\in [0,T].\end{equation}
Let
$\{P_{s,t}^{\si}\}_{T\ge t_1\ge t\ge s\ge 0}$ be the Neumann semigroup on $\bar D$ generated by
$L_t^{\si};$   that is, for any $\varphi\in C_b^2(\bar D)$,  and any  $t\in (0,T]$, $(P_{s,t}^{\si}\varphi)_{s\in [0,t]}$ is the unique solution of the PDE
\beq\label{NMM2} \pp_s u_s= - L_s^{\si} u_s,\ \ \nn_{\n}u_s|_{\pp D}=0\ \text{for}\ s\in [0,t), u_t=\varphi.\end{equation}
For any $t>0,$ let $C_b^{1,2}([0,t]\times\bar D)$ be the set of  functions  $f\in C_b([0,t]\times\bar D)$  with bounded and continuous derivatives $\pp_t f, \nn f$ and $\nn^2 f$.

 \beg{enumerate}
\item[{\bf (A4)}]    $D\in C_b^{2,L}$, {\bf (A3)} holds with   $V$  satisfying \eqref{BDR} holds for some $r_0>0$.  Moreover,
 for any  $\varphi\in C_b^2(\bar D)$ and $t\in (0,T]$, the PDE \eqref{NMM2} has a unique solution $P_{\cdot,t}^{\si}\varphi\in C_b^{1,2}([0,t]\times \bar D),$
 such that    for some  constant $c>0$  we have
 \beq\label{AB10} \|\nn^i P_{s,t}^{\si} \varphi\|_{\infty} \le  c (t-s)^{-\ff 1 2} \|\nn^{i-1} \varphi\|_{\infty},\ \ 0\le s<t\le T, \ i=1,2,\varphi\in C_b^2(\bar D).\end{equation}
\end{enumerate}

 We have the following result, where the well-posedness for $D=\R^d$ has been addressed in \cite{Ren}.

 \beg{lem} \label{L2} Assume {\bf (A3)} for $D=\R^d$ and {\bf (A4)} for $D\ne \R^d$.
 Then $\eqref{E1.1} $ is well-posed up to time $T$. Moreover, for any $t\in (0,T]$,
 \beq\label{STF} \lim_{\bar D\ni y\to x }  \|P_t^*\dd_x-P_t^*\dd_y\|_{var}=0,\ \ t \in (0,T], x\in \bar D,\end{equation}
   and $P_t$ has probability density (i.e. heat kernel) $p_t(x,y)$ such that
 \beq\label{HTK}   \inf_{x,y\in \bar D\cap B_N, \   \rr_\pp(y)\ge N^{-1}} p_t(x,y)>0,\ \ N>1, t\in (0,T], \end{equation}
 where $\inf \emptyset:=\infty$.
\end{lem}

\beg{proof}  (a)  The well-posedness.  For any $n\ge 1$, let
$$b^n:= 1_{B_n}b^{(1)}+b^{(0)}.$$
Since $b^{(1)}$ is locally bounded,
by   \cite[Theorem 1.1]{XXZZ} for $D=\R^d$ and   \cite[Theorem 2.2]{W21b} for $D\ne \R^d$,   for any $x\in \bar D$,
the following SDE is well-posed:
$$\d X_t^{x,n}= b^n(X_t^{x,n})\d t +\si(X_t^{x,n})\d W_t+\n(X_t^{x,n})\d l_t^{x,n},\ \ X_0^{x,n}=x.$$
Let $\tau_n^x:= \inf\{t\ge 0: |X_t^{x,n}|\ge n\}.$ Then $X_t^{x,n}$ solves \eqref{E1'} up to time $\tau_n^x$, and by the uniqueness we have
$$X_t^{x,n}=X_t^{x,m},\ \ t\le \tau_n^x\land\tau_m^x, n,m\ge 1.$$
So, it suffices to prove that $\tau_n^x\to\infty$ as $n\to\infty$.

Let $L_t^0:= L_t^\si+\nn_{b_t^{(0)}}.$ By \cite[Theorem 3.1]{XXZZ} for $D=\R^d$ and   \cite[Lemma 2.6]{W21b}  for $D\ne \R^d$,       {\bf (A3)} implies that  for any $\ll\ge 0$,  the PDE
\beq\label{PDE0} (\pp_t+L_t^0 )u_t=\ll u_t-b_t^{(0)},\ \ t\in [0,T], u_T=0, \nn_\n u_t|_{\pp D}=0\end{equation}
has a unique solution $u\in \tt H_{q_0}^{p_0} (T)$, and there exist constants $\ll_0, c,\theta>0$ such that
\beq\label{ESS} \ll^{\theta} (\|u\|_\infty+\|\nn u\|_\infty) +\|\pp_t u\|_{\tt L_{q_0}^{p_0}(T)}+\|\nn^2 u\|_{\tt L_{q_0}^{p_0}(T)}\le c,\ \ \ll\ge \ll_0.\end{equation}
So, we may take $\ll\ge \ll_0$ such that
\beq\label{*4}\|u\|_\infty+\|\nn u\|_\infty\le \vv,\end{equation}
where we take $\vv\le r_0$ when $\pp D$ exists.
Let $\Theta_t(x)= x+u_t(x)$.    By \eqref{BDR} and  \eqref{*4} for $\vv\le r_0$ when $\pp D$ exists,  we have
$$\<\nn V(Y_t^{x,n}), \n(X_t^{x,n})\>\d l_t^{x,n}\le 0.$$  So,  by It\^o's formula,  $Y_t^{x,n}:= \Theta_t(X_t^{x,n})$ satisfies
\beq\label{YNN0}  \d Y_t^{x,n} = \big\{1_{B_n}b_t^{(1)} +\ll u_t+1_{B_n}\nn_{b_t^{(1)}} u_t  \big\}(X_t^{x,n}) \d t + \{(\nn \Theta_t)\si_t\}(X_t^{x,n}) \d W_t+\n(X_t^n)\d l_t^n,\end{equation}
where we have used the fact that $\nn_\n u_t|_{\pp D}=0$ implies that $\{\nn \Theta_t\}\n =\n$ holds on $\pp D$. 
By \eqref{*4} and {\bf (A3)}(3) with \eqref{BDR}  when $\pp D\ne\emptyset$, there exists a constant   $c_0 >0$ such that for some martingale $M_t$,
\beg{align*} &   \d  \{V(Y_t^{x,n}) +M_t\} \\
&\le \Big[ \big\<\{b^{(1)}+\nn_{b^{(1)}} u_t\}(X_t^{x,n}), \nn  V(Y_t^{x,n})\big\>+  c_0(|\nn V(Y_t^{x,n})|+ \|\nn^2 V(Y_t^{x,n})\|)\Big]\d t \\
&\le \Big\{\<b^{(1)}(X_t^{x,n}), \nn V(X_t^{x,n})\> +   \vv |b^{(1)}(X_t^{x,n})| \sup_{B(X_t^{x,n},\vv)} \|\nn^2 V\| +c_0K V(Y_t^{x,n})\Big\}\d t\\
&\le \big\{K \phi(V(X_t^{x,n})) + c_0K V(Y_t^{x,n})\big\}\d t \le K\big\{\phi((1+\vv K) V(Y_t^{x,n})) + c_0 V(Y_t^{x,n})\big\}\d t,\ \ t\le\tau_n^x.
  \end{align*}
  Letting  $H(r):=\int_0^r\ff{\d s}{r+ \phi((1+\vv K)s)},$  by It\^o's formula and noting that $\phi'\ge 0$, we find a constant $c_1>0$ such that
 $$\d H(V(Y_t^{x,n}))\le c_1\d t +\d \tt M_t,\ \ t\in [0,\tau_n^x]$$  holds for some martingale $\tt M_t$. Thus,
 $$\E [(H\circ V)(Y_{t\land\tau_n^x}^{x,n})]\le V(x+u(x)) +c_1 t,\ \ t\ge 0, n\ge 1.$$
  Since  \eqref{*4} and  $|z|\ge n$ imply $|\Theta_t(z)|\ge |z|- |u(z)|\ge n-\vv,$  we derive
\beq\label{XXF}  \P(\tau_n^x\le t) \le  \ff{ V(x+\Theta_0(x)) +c_1 t}{\inf_{|y|\ge n-\vv} H(V(y))}=:\vv_{t,n}(x),\ \ t>0.\end{equation}
Since $\lim_{|x|\to\infty} H(V)(x)=\int_0^\infty\ff{\d s}{s+\phi((1+\vv K)s)}=\infty,$ we obtain $\tau_n^x\to\infty (n\to\infty)$ as desired.

 (b) Proof of \eqref{STF}.  By \cite[Proposition 1.3.8]{W13}, the log-Harnack inequality
 $$P_t\log f(y)\le \log P_t f(x) +c |x-y|^2,\ \ x,y\in \bar D, 0<f\in \B_b(\bar D)$$ for some constant $c>0$   implies the gradient estimate
 $$|\nn P_t f|^2 \le 2 c P_t |f|^2,\ \ f\in \B_b(\bar D),$$ and hence
 $$\lim_{y\to x} \|P_t^*\dd_x-P_t^*\dd_y\|_{var}=0,\ \ x\in \bar D.$$
 Let $P_t^{n}$ be the Markov semigroup associated with $X_t^{n}$. Thus, by
the log-Harnack inequality in \cite[Theorem 4.1]{YZ0} for $D=\R^d$ and in  \cite[Theorem 4.1]{W21b} for $D\ne \R^d$,  we have
 \beq\label{STF'} \lim_{y\to x} \|(P_t^n)^*\dd_x-(P_t^n)^*\dd_y\|_{var}=0,\ \ t\in (0,T].\end{equation}
On the other hand, by \eqref{XXF} and $X_t=X_t^n$ for $t\le \tau_n$, we obtain
\beg{align*}&\lim_{n\to\infty} \sup_{y\in \bar D\cap B(x,1)} \|P_t^*\dd_y-(P_t^n)^*\dd_y\|_{var}=\lim_{n\to\infty} \sup_{|f|\le 1,y\in \bar D\cap B(x,1)} |P_tf(y)-P_t^n f(y)|\\
&\le 2 \lim_{n\to\infty} \sup_{y\in \bar D\cap B(x,1)}\P(\tau_n^y\le t)=0.\end{align*}
Combining this with \eqref{STF'} and the triangle inequality, we prove \eqref{STF}.

 (c) Finally,  let $L_t:=L_t^\si+\nn_{b_t}$. By It\^o's formula, for any $f\in C_0^2((0,T)\times D)$ we have
 $$\d f_t(X_t)= (\pp_t+L_t)f_t(X_t)\d t+ \d M_t$$ for some martingale $M_t$, so that $f_0=f_T=0$ yields
 $$\int_{(0,T)} P_t \{(\pp_t+L)f_t\}\d t  =0,\ \ f\in C_0^\infty((0,T)\times D).$$This implies that the heat kernel $p_t(x,\cdot)$ of $P_t$ solves the following PDE
 on $(0,T)\times D$ in the weak sense: 
 $$\pp_t u_t = L_t^* u_t= {\rm div} \scr A(t,\cdot, u_t,\nn u_t) +\scr B(t,\cdot,\nn u_t),$$
 where $\scr A:=(\scr A_1,\cdots, \scr A_d)$ and $\scr B$ are defined as
\beg{align*} &\scr A_i(t,\cdot,u,\nn u):=  \ff 1 2 \sum_{j=1}^d (\si_t\si_t^*)_{ij}  \pp_j u +\sum_{j=1}^d \Big\{\ff  1 2\pp_j (\si_t\si_t^*)_{ij} -b_t^i\Big\} u,\\
&\scr B(t,\cdot, \nn u):= - \sum_{i,j=1}^d \big\{\pp_j (\si_t\si_t^*)_{ij}\big\}\pp_i u.\end{align*}
 By   the Harnack inequality   as in \cite[Theorem 3]{Serr} (see also \cite{Trud}), under the given conditions,  
 for any $0<s<t\le T$ and $N>1$ with
 $$\tt B_N:= \big\{x\in \bar D \cap B_N: \ \rr_\pp(x) \ge N^{-1}\big\}$$ having positive volume,
 there exists a constant $c(s,t,N)>0$ such that  satisfies
 \beq\label{HNK} \sup_{\tt B_N} p_s(x,\cdot)\le c(s,t,N) \inf_{\tt B_N} p_{t} (x,\cdot),\ \ x\in \bar D.\end{equation}
Since $\int_{\tt B_N} p_s(x,y)\d y\to 1$ as $N\to\infty$, this  implies $p_t(x,y)>0$ for any $(t,x,y)\in (0,T]\times \bar D\times D$. In particular,
$P_t 1_{\tt B_N}>0$. On the other hand, \eqref{STF} implies that $P_t 1_{\tt B_N}$ is continuous, so that
  $$\inf_{x\in \bar D\cap  B_N} P_t 1_{\tt B_N}(x)>0,\ \ t\in (0,T].$$   This together with \eqref{HNK}
 gives
 $$\inf_{(\bar D\cap  B_N)\times \tt B_N} p_{t} \ge\ff 1 { c(s,t,N)  }  \inf_{x\in \bar D\cap \bar B_N} P_s 1_{\tt B_N}(x)  >0,\ \ 0<s<t\le T.$$
Therefore,   \eqref{HTK} holds.

\end{proof}

To make Zvonkin's transform to kill the singular drift, we present the  lemma which extends  Theorem 2.10 in \cite{XZ} for $D=\R^d$.

\beg{enumerate}
\item[{\bf (A5)}]    $D=\R^d$,   $\si$ and $b^{(0)}$ satisfy the following conditions.
\item[$(1)$]     $a:= \si\si^*$ is invertible and uniformly continuous with $\|a\|_\infty+\|a^{-1}\|_\infty<\infty$.
\item[$(2)$]      $|b^{(0)}|\in \tt L^{p}$ for some  $p>d$.
 \end{enumerate}

 \beg{enumerate}
\item[{\bf (A6)}]    $\pp D\in C_b^{2,L}$,   {\bf (A5)} holds for $\bar D$ replacing $\R^d$, and {\bf (A2)}(2) holds.
 \end{enumerate}

\beg{lem} \label{L1} Assume {\bf (A5)} for $D=\R^d$ and {\bf (A6)} for $D\ne \R^d$. Let $L^0=\ff 1 2 {\rm tr}\{\si\si^*\nn^2\} +\nn_{b^{(0)}}.$ Then there exist constants  $\ll_0>0$ increasing in $\|b^{(0)}\|_{\tt L^{p}}$ such that
for any $\ll\ge \ll_0$ and any $f\in \tt L^{k}$ for some $k \in (1,\infty)$,
the elliptic equation
\beq\label{*1} (L^0-\ll)u= f,\ \nn_\n u|_{\pp D}=0\ \text{if}\ D\ne \R^d \end{equation}
has a unique solution $u\in \tt H^{2,k}$. Moreover, for any $p'\in [k,\infty]$ and $\theta \in [0, 2-\ff d {k}+\ff d{p'})$,   there exists a  constant $c>0$ increasing in $\|b^{(0)}\|_{\tt L^{p}}$ such that
\beq\label{*2} \ll^{\ff 1 2 (2-\theta +\ff d {p'}-\ff d {k})}\|u\|_{\tt H^{\theta, p'}} +\|u\|_{\tt H^{2,k}}\le c\|f\|_{\tt L^{k}},\ \ f\in\tt L^{k}.\end{equation}
\end{lem}

\beg{proof}  (a) Let us verify the priori estimate \eqref{*2} for a solution $u$ to \eqref{*1}, which in particular implies the uniqueness, since the difference of two solutions solves the equation with $f=0$.

For $u\in \tt H^{2,k}$ solving \eqref{*1},  let
$$\bar u_t= u (1-t),\ \ t\in [0,1].$$ By \eqref{*1} we have
$$(\pp_t +L^0 -\ll)\bar u_t= f(1-t)- u,\ \ t\in [0,1], \bar u_1=0,\ \nn_\n \bar u_t|_{\pp D}=0\ \text{if}\ D\ne \R^d.$$
By Theorem 2.1 with $q=q'=2$ in \cite{YZ0} for $D=\R^d$, and Lemma 2.6 in \cite{W21b} for $D\ne \R^d$,
there exist  constants $\ll_1, c_1>1$ increasing in $\|b^{(0)}\|_{\tt L^{p}}$ and sufficient large $q>2$ such that
\beq\label{*3} \ll^{\ff 1 2 (2-\theta +\ff d {p'}-\ff d k)}\|\bar u\|_{\tt H^{\theta, p'}_q} +\|\bar u\|_{\tt H^{2,k}_q}\le c_1\|f(1-t) -u\|_{\tt L^k_q}\le c_1 \|f\|_{\tt L^k} +   c_1\|u\|_{\tt L^k}.\end{equation}
Taking $\theta=0, p=p'$ and $c_2= \|1-\cdot\|_{L^q([0,1])}$, we obtain
$$ \ll^{\ff 1 2 (2-\theta ) }\|u\|_{\tt L^p}\le \ff{ c_1}{c_2} \big(\|f\|_{\tt L^p} +   \|u\|_{\tt L^p}\big),\ \ \ll\ge \ll_1.$$
Letting    $\ll_0>\ll_1$ such that
$$ \ll_0^{\ff 1 2 (2-\theta)}\ge 2 \ff{c_1}{c_2},$$ we obtain we obtain
$$\|u\|_{\tt L^k}\le \|f\|_{\tt L^k},\ \ \ll\ge \ll_0.$$  Combining this with  \eqref{*3} implies \eqref{*2} for some constant $c>0$.

(b) Existence of solution for $f\in \tt L^k$. Let $\{f_n\}_{n\ge 1}\subset C_b^\infty (\bar D)$ with $\n f_n|_{\pp D}=0$ is $\pp D\ne \emptyset$ such that $\|f_n-f\|_{\tt L^k}\to 0$ as $n\to\infty$. Let $P_t^0$ be the Markov semigroup generated by $L^0$ and  let
$$u_n= \int_0^\infty \e^{-\ll t}P_t^0 f_n \d t.$$ By Kolmogorov equation we have
$$\pp_t P_t^0f_n = L^0P_t^0 f_n= P_t^0 L^0 f_n$$ so that
$$L^0 u_n =\int_0^\infty \e^{-\ll t} L^0P_t^0f_n\d t= \int_0^\infty \e^{-\ll t} \pp_t P_t^0 f_n\d t= \ll u_n-f_n.$$
Then
$$(L^0-\ll)(u_n-u_m)= f_n-f_m,\ \ n,m\ge 1.$$
By \eqref{*2},
$$\lim_{n,m\to\infty} \big\{\|u_n-u_m\|_{\tt H^{\theta, p'}}  +\|\nn^2(u_n-u_m)\|_{\tt L^k} \big\}=0,$$
so that $u:=\lim_{n\to\infty} u_n$ exists in $\tt H^{\theta,p'}\cap \tt H^{2,k}$, which solves \eqref{*1}.
\end{proof}

  \subsection{Proofs of Theorems \ref{T1.6.1}, \ref{T2} and Corollary \ref{C1.2}}
  \beg{proof}[Proofs of Theorems \ref{T1.6.1}, \ref{T2}]
It is easy to see that  \eqref{LYP} implies {\bf (A3)}(3) for any $T>0$ and $\phi(r)=1$,    by Lemma \ref{L2}, {\bf (A1)} and {\bf (A2)} imply the well-posedness, strong Feller property and irreducibility of  \eqref{E1'} and \eqref{E1'b} respectively.
 According to \cite[Theorem 4.2.1]{DZ96},  the strong Feller property and the irreducibility   imply the uniqueness of invariant probability measure. So, it remains to prove the existence  of the invariant probability measure $\mu$ and the claimed  assertions on  the ergodicity.

(a)  Let $u$ solve \eqref{*1} for $b=-b^{(0)} $ and large enough $\ll>0$ such that \eqref{*2} implies \eqref{*4}. Moreover, for
  $\Theta(x):=x+u(x)$, let  $\hat P_t$ be the Markov semigroup associated with $Y_t:= \Theta(X_t)$, so that
\beq \label{HAT} \hat P_tf(x)= \{P_t(f\circ \Theta)\}(\Theta^{-1}(x)),\ \ t\ge 0, x\in\R^d, f\in \B_b(\R^d).\end{equation}
Since $\lim_{|x|\to\infty} \sup_{|y-x|\le \vv} \ff{|\nn V(y)|}{V(x)}=0,$ by \eqref{*4} and $V\ge 1$ we find a constant $\theta\in (0,1)$  such that
\beq\label{*GMM}\|\nn u(x)\| \lor  |\Theta(x)-x|\le \vv,\ \ \theta V(\Theta(x))\le V(x)\le \theta^{-1} V(\Theta(x)),\ \ \   x\in\bar D.\end{equation}
Thus, it  suffices to prove the desired assertions for $\hat P_t$ replacing $P_t$, where the unique invariant probability measure   $\hat\mu$ of $\hat P_t$ and that $\mu$ of $P_t$
satisfies
\beq\label{INV} \hat\mu=\mu\circ\Theta^{-1}. \end{equation}

(b) Let $X_t^n, Y_t^n$ and $\tau_n$ be in the proof of Lemma \ref{L2} for the present time-homogenous setting.
Since $Y_t^n=Y_t$ and $ 1_{B_n}(X_t^n)=1$ for  $t\le \tau_n$,  and since $\tau_n\to\infty$ as $n\to\infty$, \eqref{YNN0} implies
$$ \d Y_t = \big\{b^{(1)} +\ll u+\nn_{b^{(1)}} u  \big\}(X_t) \d t + \{(\nn \Theta)\si\}(X_t) \d W_t+\n(X_t)\d l_t,$$
so that  for any $\vv\in (0,1\land r_0)$, where $r_0>0$ is in \eqref{BDR} when $\pp D\ne \emptyset$, by It\^o's formula  and \eqref{BDR}, we find a constant $c_\vv>0$  such that
\beg{align*} &   \d  \{V(Y_t) +M_t\} \le \Big\{\big\<\{b^{(1)}+\nn_{b^{(1)}} u\}(X_t), \nn  V(Y_t)\big\>+  c_\vv(|\nn V(Y_t)|+ \|\nn^2 V(Y_t)\|)\Big\}\d t \\
&\le \Big\{\<b^{(1)}(X_t), \nn V(X_t)\> +   \vv |b^{(1)}(X_t)| \sup_{B(X_t, \vv)}\{|\nn V|+\|\nn^2 V\| \}+  c_\vv\sup_{B(X_t, \vv)} (|\nn V|+ \|\nn^2 V\|)\Big\}\d t.  \end{align*}
Combining this with \eqref{LYP} and \eqref{BDR}  for $D\ne \R^d$, when $\vv>0$ is small enough we find   constants $c_1,c_2>0$ such that
$$ \d  \{V(Y_t) +M_t\} \le \{c_1 -c_2\Phi(V(X_t))\}\d t.$$
By \eqref{*GMM},  this implies that for some constant   $c_4>0,$
\beq\label{GMM} \d V(Y_t)\le \big\{c_4-c_2\Phi(\theta V(Y_t))\big\}\d t-\d M_t.\end{equation}
Thus,
$$\int_0^t \E \Phi(\theta V(Y_s))\d s\le \ff {c_4+V(x)}{c_2}<\infty,\ \ t>0, Y_0=x\in \Theta(\bar D).$$
Since $\Phi(\theta V)$ is a compact function, this implies the existence of invariant probability $\hat \mu$ according to the standard Bogoliubov-Krylov's   tightness argument.
Moreover,   \eqref{GMM} implies $\hat\mu(\Phi(\theta V))<\infty$, so that by \eqref{*GMM} and \eqref{INV}, $\mu(\Phi(\vv_0 V))<\infty$ holds for $\vv_0=\theta^2.$

 (c) By  \eqref{HTK}, \eqref{HAT} and \eqref{*GMM},    any compact set $\mathbf K\subset \Theta(\bar D) $ is a petite set of $\hat P_t$, i.e. there exist $t>0$ and a nontrivial measure $\nu$ such that
 $$ \inf_{x\in \mathbf K} \hat P_t^*\dd_x\ge \nu. $$
   When $\Phi(r)\ge kr$ for some constant $k>0$,     \eqref{GMM} implies
 \beq\label{ACC} \hat P_t  V(x)\le \ff{k_1}{k_2}+ \e^{-k_2 t}   V(x),\ \ t\ge 0, x\in \Theta(\bar D)\end{equation}
 for some constants $k_1,k_2>0$. Since $\lim_{|x|\to\infty}V(x)=\infty$ and as observed above that any compact set is a petite set for $\hat P_t$,
 by Theorem 5.2(c) in \cite{MT}, we obtain
 $$  \|\hat P_t^*\dd_x-\hat\mu\|_{V}  \le c\e^{-\ll t} V(x),\ \ x\in \Theta(\bar D), t\ge 0$$
 for some constants $c,\ll>0$.  Thus,
 $$ \|\hat P_t^*\dd_x-\hat P_t^*\dd_y\|_{V}  \le c\e^{-\ll t} (V(x)+V(y)),\ \ t\ge 0, x,y\in \Theta(\bar D).$$
 Therefore, for any probability measures $\mu_1,\mu_2$ on $\Theta(\bar D)$,
 \beg{align*}&\|\hat P_t^*\mu_1-\hat P_t^*\mu_2\|_{V} = \|\hat P_t^*(\mu_1-\mu_2)^+ - \hat P_t^*(\mu_1-\mu_2)^-\|_V\\
 &=\ff 1 2 \|\mu_1-\mu_2\|_{var} \Big\|\hat P_t^*\ff{2(\mu_1-\mu_2)^+}{\|\mu_1-\mu_2\|_{var}}  - \hat P_t^*\ff{2(\mu_1-\mu_2)^-}{\|\mu_1-\mu_2\|_{var}}\Big\|_V\\
 &\le \ff c 2\e^{-\ll t} \|\mu_1-\mu_2\|_{var}  \Big(\ff{2(\mu_1-\mu_2)^+}{\|\mu_1-\mu_2\|_{var}}  + \ff{2(\mu_1-\mu_2)^-}{\|\mu_1-\mu_2\|_{var}}\Big)(V)\\
 &\le  c \e^{-\ll t} \|\mu_1-\mu_2\|_V.\end{align*} This together with  \eqref{HAT} and \eqref{*GMM}  implies
 \eqref{EX1} for some constants $c,\ll>0$.

 (d) Let  $\Phi$ be convex.  By Jensen's inequality and \eqref{GMM},
$\gg_t:=  \theta \E [V(Y_t)]$   satisfies
\beq\label{GGM}  \ff{\d}{\d t} \gg_t \le \theta c_4- \theta  c_2 \Phi( \gg_t),\ \ t\ge 0.\end{equation}
Let
$$H(r):=\int_0^r \ff{\d s}{\Phi(s)},\ \ r\ge 0.$$
We aim to prove that for some constant $k>1$
\beq\label{GMM2} \gg_t\le   k+H^{-1}(H(\gg_0)-tk^{-1}),\ \ \ t\ge 0,\end{equation}
where $H^{-1}(r):=0$ for $r\le 0$.
We prove this estimate by considering three situations.
\beg{enumerate}
\item[(1)] Let  $\Phi(\gg_0)\le \ff{c_4}{c_2}.$ Since  \eqref{GGM} implies  $\gg_t'\le 0$ for $ \gg_t \ge \Phi^{-1} ( \ff{c_4}{c_2})$, so
\beq\label{FMM3}  \gg_t \le \Phi^{-1} ( c_4/c_2),\ \ t\ge 0.\end{equation}
\item[(2)]  Let  $\ff{c_4}{c_2}<\Phi(\gg_0)\le \ff{2c_4}{c_2}.$  Then  \eqref{GGM} implies $\gg_t'\le 0$ for all $t\ge 0$ so that
\beq\label{FMM4}  \gg_t \le \Phi^{-1} ( 2c_4/c_2),\ \ t\ge 0.\end{equation}
\item[(3)] Let $\Phi(\gg_0)> \ff{2c_4}{c_2}.$ If
$$t\le t_0:=\inf\Big\{t\ge 0: \Phi(\gg_t)\le  \ff{2c_4}{c_2}\Big\},$$
  then \eqref{GGM} implies
$$\ff{\d H(\gg_t)}{\d t}=\ff{\gg_t'}{\Phi(\gg_t)}\le -\ff{ \theta  c_2}2 , $$ so that
\beq\label{TXV}  H(\gg_t)\le H(\gg_0)- \ff{ \theta  c_2}2 t,\ \ t\in [0, t_0],\end{equation}
which implies
$$\gg_t\le H^{-1}(H(\gg_0)- \theta c_2 t/2),\ \ t\in [0,t_0].$$
Noting that when $t>t_0$, $(\gg_t)_{t\ge t_0} $ satisfies \eqref{GGM} with $\gg_{t_0}$ satisfies $\ff{c_4}{c_2}<\Phi(\gg_{t_0})\le \ff{2c_4}{c_2},$
so that \eqref{FMM3} holds, i.e.
$$\gg_t\le \Phi^{-1} ( 2c_4/c_2).$$ In conclusion, we obtain
$$\gg_t\le \Phi^{-1} ( 2c_4/c_2)+ H^{-1}(H(\gg_0)- \theta c_2t/2),\ \ t\ge 0.$$
Combining this with (1) and (2), we prove \eqref{GMM2} for some constant $k>1$. \end{enumerate}

(e) Since $1\le \Phi(r)\to\infty$ as $r\to\infty$, when $\Phi$ is convex we find a constant $\dd>0$ such that $\Phi(r)\ge \dd r, r\ge 0$. So, by step (b), \eqref{EX1} holds.
Combining this with \eqref{GMM2} and applying the semigroup property, we derive
\beg{align*}& \|\hat P_{t}^*\dd_x-\hat\mu\|_{V}=\sup_{|f|\le V} |\hat P_{t/2}(\hat P_{t /2}f-\hat\mu(f))(x)|\\
& \le c \e^{-\ll t/2} \hat P_{t/2} V(x) \le c\big\{k+H^{-1}(H(\theta V(x))-(2k)^{-1} t)\big\}\e^{-\ll t/2}.\end{align*}
Combining this with \eqref{HAT}, \eqref{*GMM} and \eqref{INV}, we prove
  \eqref{EX0}   for some constants $k,\ll>0.$

Finally, if $H(\infty)<\infty$, we take $t^*= k H(\infty)$ in \eqref{EX0} to derive
$$\sup_{x\in \bar D} \|P_t\dd_x- \mu\|_V\le c \e^{-\ll t},\ \ t\ge t^*$$ for some constants $c,\ll>0$, which
implies  \eqref{EX2} by the argument  leading to \eqref{EX1} in step (c).  \end{proof}

\beg{proof}[ Proof of Corollary \ref{C1.2}]
 By \eqref{PII}, for any $\theta\in ((1-\aa)^+,\ff 1 2)$ there exists a constant $c_3>0$ such that
$$\phi(r)\ge c_3 (1+r)^{1-\theta},\ \ r\ge 0.$$
Then \eqref{LYP} holds for $V:= \e^{(1+|\cdot|^2)^\theta}$ and $\Phi(r)=r$. So the first assertion in (1) follows from Theorem \ref{T1.6.1}(1).

Next, \eqref{B1} and \eqref{PSI} imply \eqref{LYP} for $V:= \psi(|\cdot|^2)$ and $\Phi(r)=r$, so that the second assertion in (1) holds by  Theorem \ref{T1.6.1}(1).

Finally, if $\int_0^\infty\ff{\d s}{\phi(s)}<\infty$, then for any $q>0$, \eqref{LYP} holds for $V:=(1+|\cdot|^2)^q$ and $\Phi(r)= (1+r)^{1-\ff 1 q}\phi(r^{\ff 1 q})$, so that
$\int_0^\infty\ff{\d s}{\Phi(s)}<\infty$. Then the proof is finished by Theorem \ref{T1.6.1}(2). \end{proof}

 \section{Uniform ergodicity for singular reflecting  McKean-Vlasov SDEs  }

 We now consider the SDE \eqref{E1} for $D=\R^d$ or $D$ being a $C_b^{2,L}$ domain. To prove the uniform ergodicity, compare \eqref{E1} with the following classical SDE with fixed distribution parameter  $\gg\in \scr P$:
\beq\label{EGG} \d X_t^\gg= b(X_t^\gg,\gg)+ \si(X_t^\gg) \d W_t+\n(X_t^\gg)\d l_t^\gg,\end{equation}
where we set $l_t^\gg=0$ if $D=\R^d.$ 
If for any $\nu\in \scr P$,  {\bf (A1)} for $D=\R^d$ or {\bf (A2)} for $D\ne \R^d$ holds for $b(\cdot,\nu)$ replacing $b$, then  the well-posedness follows from that of \eqref{EGG} and  \cite[Theorem 3.2]{W21b} for $k=0$.
Let $(P_t^\gg)^*\nu=\L_{X_t^\gg}$ for $\L_{X_0^\gg}=\nu$. 

Let $\zeta(\gg_1,\gg_2):= \si^*(\si\si^*)^{-1}\big[b(\cdot,\gg_2)-b(\cdot,\gg_1)\big], \ \ \gg_1,\gg_2\in\scr P.$  We make the following assumption on the dependence of distribution.

 \beg{enumerate} \item[{\bf (H1)}]  {\bf (A1)} for $D=\R^d$ or {\bf (A2)} for $D\ne \R^d$ holds with  $b(\cdot,\nu)$ replacing $b$  uniformly in $\nu\in\scr P$. 
 \item[{\bf (H2)}] There exist constants $q\ge 2$ and $k>0$
 such that  for any $\gg\in \scr P$ and $\nu\in C([0,1]; \scr P)$,
\beg{align*}& \int_0^t \d s \int_{\bar D} |\zeta(\gg,\nu_s) |^2\d (P_s^\gg)^*\nu_0\le k^2 \bigg(\int_0^t \|\nu_s-\gg\|_{var}^q\d s\bigg)^{\ff 2 q},\\
&\int_{\bar D} \e^{\ff 1 2 \int_0^t  |\zeta(\gg,\nu_s) |^2}\d (P_s^\gg)^*\nu_0<\infty,\ \ \ t\ge 0.\end{align*} 
\end{enumerate} 
\paragraph{Remark 3.1.} (1) Obviously, {\bf (H2)}  holds for $q=2$ if there exists a constant $\kk>0$ such that 
$$ |b(x,\gg_1)-b(x,\gg_2)|\le \kk \|\gg_1-\gg_2\|_{var},\ \ x\in\bar D, \gg_1,\gg_2\in \scr P.$$
 By a standard fixed point argument, this together with $(H_1)$ implies  the well-posedness of  $\eqref{E1}$ for any initial value, see \cite{W21b}.    

(2)  In general, {\bf (H2)} follows from Krylov's estimate for $p,q>2$ with $\ff d p+\ff 2 q<1$ and the condition
 $$\|\zeta(\gg_1,\gg_2)\|_{\tt L^p}\le k \|\gg_1-\gg_2\|_{var},\ \ \ \gg_1,\gg_2\in\scr P.$$
 Note that under {\bf (H1)},  the Krylov's estimate  holds when $b$ contains an $\tt L^p$ term for $p>d$ and a Lipschitz continuous term, see \cite{W21b} for details. 

 \beg{thm}\label{T3} Assume  {\bf (H1)} and {\bf (H2)} and let   $\eqref{E1}$ be well-posed. 
  If $k$ is small enough and  $\Phi$ is convex with $\int_0^\infty \ff{\d s}{\Phi(s)}<\infty$, then
 $P_t^*$ has a unique invariant probability measure $\mu$,  $\mu(\Phi(\vv_0 V))<\infty$ holds for some constant $\vv_0>0$,  and  there exist constants $c,\ll>0$ such that
 \beq\label{Q2} \|P_t^*\nu -\mu \|_{var}\le c \e^{-\ll t}\|\mu-\nu\|_{var},\ \ t\ge 0, \nu\in \scr P.\end{equation}  
  \end{thm}

To prove this result, we first present a general result deducing the uniform ergodicity of McKean-Vlasov SDEs from that of classical ones.

  The following result says that if \eqref{EGG} is uniformly  ergodic uniformly in $\gg$, and if the dependence of $b(x,\mu)$ on $\mu$ is weak enough, then
\eqref{E1} is uniformly ergodic.  

\beg{lem}\label{L3} Assume   $(H_1)$ and that for any $\gg\in \scr P$,    $(P_t^\gg)^*$ has a unique invariant probability measure
$\mu_\gg$ such that
\beq\label{LNA} \|(P_t^\gg)^*\mu -  \mu_\gg  \|_{var} \le c \e^{-\ll t} \|\mu -\mu_\gg\|_{var},\ \ t\ge 0, \gg, \mu  \in \scr P\end{equation}
holds for some constants $c,\ll>0$.  
Then $\eqref{E1}$ is well-posed and the following assertions hold.
\beg{enumerate} \item[$(1)$]  If there exists a constant $\kk\in (0,\kk_1)$ for
$$ \kk_1:=\sup_{t>(\log c)/\ll} \ff{1-c\e^{-\ll t}}{\ss t},$$  such that 
\beq\label{DST1} \int_{\bar D} |\zeta(\gg_1,\gg_2)|^2 \d\mu_{\gg_2} \le \kk^2\|\gg_1-\gg_2\|_{var}^2,\ \ \gg_1,\gg_2\in \scr P,\end{equation}  
then $P_t^*$ associated with $\eqref{E1}$ has a unique invariant probability measure $\mu$.
\item[$(2)$] Let $\mu$ be $P_t^*$-invariant. If  there exist constants $q\ge 2$ and  $k\in (0,k_q)$, where 
$$ k_q:= \sup\bigg\{k>0:\ \ff {4^{q-1}(ck)^q \e^{2^{q-1}k^q t} }{q\ll+ 2^{q-1}k^q}\le \ff 1 2\bigg\},$$
such that for $\hat t:=\ff{\log (2c)}\ll$, 
\beq\label{DST2}\E \int_0^t   |\zeta(\mu, P_s^*\nu)(X_s^\mu)|^2  \d s\le k^2 \bigg(\int_0^t \|\mu-P_s^*\nu\|_{var}^q\d s\bigg)^{\ff 1 q},\ \ t\in (0,\hat t], \nu\in \scr P,\end{equation} 
then   there exists a  constant $c' >0$ such that
\beq\label{Q'} \|P_t^*\nu-\mu\|_{var}^q\le c'\e^{-\ll' t}\|\nu-\mu\|_{var}^q,\ \ t\ge 0, \nu\in \scr P\end{equation}
holds for
$$\ll':= -\ff{\ll}{\log (2c)} \log\Big(\ff 1 2+\ff {4^{q-1}(ck)^q \e^{2^{q-1}k^q t} }{q\ll+ 2^{q-1}k^q}\Big)>0.$$
\end{enumerate} \end{lem}

\beg{proof} 

(a) Existence and uniqueness of  $\mu$. For any $\gg\in\scr P$,
\eqref{LNA} implies that $P_t^\gg$ has a  unique invariant probability measure $\mu_\gg$. It suffices to prove that the map
$\gg\mapsto \mu_\gg$ has a unique fixed point $\mu$, which is the unique invariant probability measure of $P_t^*$.

For $\gg_1,\gg_2\in \scr P$, \eqref{EGG} implies
\beq\label{PA1} \|(P_t^{\gg_1})^*\mu_{\gg_2}- \mu_{\gg_1}\|_{var} \le c\e^{-\ll t} \|\mu_{\gg_2}-\mu_{\gg_1}\|_{var},\ \ t\ge 0.\end{equation}
On the other hand, let $(X_t^1,X_t^{2}) $ solve the SDEs
$$ \d X_t^i= b(X_t^i,\gg_i)+ \si(X_t^i) \d W_t+\n(X_t^i)\d l_t^i,\ \ i=1,2$$
with $X_0^1=X_0^2$ having distribution $\mu_{\gg_2}$. Since $\mu_{\gg_2}$ is $(P_t^{\gg_2})^*$-invariant,  we have
\beq\label{PA2} \L_{X_t^2}= (P_t^{\gg_2})^*\mu_{\gg_2}=\mu_{\gg_2},\ \
\L_{X_t^1} =(P_t^{\gg_1})^*\mu_{\gg_2},\ \ t\ge 0.\end{equation}
By $(H_2)$, 
$$R_t=\e^{\int_0^t \<\zeta(\gg_1,\gg_2)(X_s^1), \d W_s\>-\ff 1 2 \int_0^t |\zeta(\gg_1,\gg_2) (X_s^1)|^2\d s},\ \ t\ge 0$$
is a martingale, and by Girsanov's theorem, for any $t>0$,
$$\tt W_r:= W_r- \int_0^r\zeta(\gg_1,\gg_2)(X_s^1)\d s,\ \ r\in [0,t] $$
is a Brownian motion under $\Q_t:=R_t\P$. Reformulating the SDE for $X_r^1$ as
$$\d X_r^1= b(X_r^1, \gg_2)\d r+ \si(X_r^1)\d \tt W_r+\n(X_r^1)\d l_r^1,\ \ r\in [0,t],$$
by $X_0^1=X_0^2$ and the weak uniqueness, the law of $X_t^1$ under $\Q_t$ satisfies
$$\L_{X_t^1|\Q_t}= \L_{X_t^2} =(P_t^{\gg_2})^*\mu_{\gg_2}=\mu_{\gg_2}.$$
 Combining this with   \eqref{PA2} and Pinsker's inequality,  we obtain
\beq\label{PSK} \beg{split}& \|(P_t^{\gg_1})^*\mu_{\gg_2}  - \mu_{\gg_2} \|_{var}^2=\|(P_t^{\gg_1})^* \mu_{\gg_2}- (P_{t}^{\gg_2})^*\mu_{\gg_2} \|_{var}^2\\
&= \sup_{|f|\le 1} \big|\E[f(X_t^1)]- \E[f(X_t^1)R_t]\big|^2 \le \big(\E |R_t-1|\big)^2 \le 2 \E[R_t\log R_t]\\
&= 2 \E_{\Q_t} [\log R_t]=\E_{\Q_t} \int_0^t\big| \zeta(\gg_1,\gg_2)|^2(X_s^1)\big|^2\d s=t\int_{\bar D} |\zeta(\gg_1,\gg_2)|^2 \d\mu_{\gg_2}.\end{split}\end{equation}
Then \eqref{DST1} implies
$$\|(P_t^{\gg_1})^*\mu_{\gg_2}  - \mu_{\gg_2} \|_{var}^2\le \kk^2 t \|\gg_1-\gg_2\|_{var}^2.$$
Combining this with \eqref{PA1} and taking $t=\ff{\log(2c)}{\ll},$ we derive
\beg{align*} &\|\mu_{\gg_1} - \mu_{\gg_2} \|_{var}\le \|(P_t^{\gg_1})^* \mu_{\gg_2}-\mu_{\gg_1}   \|_{var}+ \|(P_t^{\gg_1})^* \mu_{\gg_2}- \mu_{\gg_2} \|_{var}\\
&\le \  c\e^{-\ll t}  \|\mu_{\gg_1}-\mu_{\gg_2}\|_{var}+ \kk\ss t \|\gg_1-\gg_2\|_{var},\ \ t>0.\end{align*}
Thus,
$$ \|\mu_{\gg_1} - \mu_{\gg_2} \|_{var}\le \inf_{t>(\log c)/\ll}  \ff{\kk \ss t}{1-c\e^{-\ll t}}  \|\gg_1-\gg_2\|_{var}\\
 =\ff{\kk}{\kk_1} \|\gg_1-\gg_2\|_{var}.$$
Since   $\kk<\kk_1,$   $\mu_\gg$ is contractive in $\gg$, hence  has a unique fixed point.

(b) Uniform ergodicity. Let $\mu$ be the unique invariant probability measure of $P_t^*$, and for any $\nu\in\scr P$ let $(\bar X_0,X_0)$ be $\F_0$-measurable such that
$$\P(\bar X_0\ne X_0)= \ff 1 2 \|\mu-\nu\|_{var},\ \ \L_{\bar X_0}=\mu, \ \ \L_{X_0}=\nu.$$
 Let $\bar X_t$ and $X_t$ solve the following SDEs with initial values $ \bar X_0$ and $X_0$ respectively:
\beg{align*} &\d\bar X_t=b(\bar X_t,\mu)\d t+\si(\bar X_t)\d W_t+\n(\bar X_t)\d\bar l_t,\\
& \d X_t=b(X_t, P_t^*\nu)\d t+\si(X_t)\d W_t+\n(X_t)\d l_t.\end{align*}
  Since $\mu$ is $P_t^*$-invariant, we have
\beq\label{UMM} \L_{\bar X_t}=(P_t^\mu)^*\mu=P_t^*\mu=\mu.\end{equation}  Moreover, $\L_{X_t}= P_t^*\nu$ by the definition of $P_t^*$.
Let
$$\bar R_t:= \e^{\int_0^t \<\zeta(P_s^*\nu,\mu)(X_s), \d W_s\>-\ff 1 2 \int_0^t |\zeta(P_s^*\nu,\mu)(X_s)|^2\d s}.$$
Similarly to \eqref{PSK}, by  Girsanov's theorem we have $\L_{X_t|\bar R_t\P}=(P_t^\mu)^*\nu$, so that 
 Pinsker's inequality and \eqref{DST2} yield 
\beg{align*} &\| (P_t^\mu)^* \nu- P_t^*\nu\|_{var}^2= \sup_{|f|\le 1} \big|\E[f(X_t)\bar R_t]-\E[f(X_t)]\big|^2 \\
&\le k^2  \bigg(\int_0^t \|\mu- P_s^*\nu\|_{var}^q\bigg)^{\ff 2 q}\d s,
\ \ t\in [0,\hat t].\end{align*} 
This together with \eqref{PA1} for $\gg_1=\mu$ and \eqref{UMM} gives
\beg{align*} &\|P_t^*\nu-\mu\|_{var}^q\le2^{q-1}\big( \|P_t^*\nu-(P_t^\mu)^*\nu\|_{var}^q+ \| (P_t^\mu)^*\nu-\mu \|_{var}^q\big)\\
&\le 2^{q-1} k^q  \int_0^t \|\mu- P_s^*\nu\|_{var}^q\d s+ 2^{q-1} c^q\e^{-q\ll t} \|\nu-\mu\|_{var}^q, \ \  t\in [0,\hat t].\end{align*}
By Gronwall's inequality we obtain
\beg{align*} &\|P_t^*\nu-\mu\|_{var}^q\le  \|\mu-\nu\|_{var}^q\bigg(2^{q-1}c^q\e^{-q\ll t}+4^{q-1} k^q c^q \int_0^t \e^{-q\ll s+2^{q-1}k^q (t-s)} \d s\bigg)\\
&\le \Big\{ 2^{q-1}c^q\e^{-q\ll t} + \ff {4^{q-1}(ck)^q \e^{2^{q-1}k^q t} }{q\ll+ 2^{q-1}k^q} \Big\} \|\mu-\nu\|_{var}^q,\ \  t\in [0,\hat t].\end{align*}
Taking $t=\hat t:=\ff{\log (2c)}\ll,$ we arrive at
$$\|P_{\hat t}^*\nu-\mu\|_{var}^2\le  \dd_k  \|\mu-\nu\|_{var}^2,\ \ \nu\in \scr P$$ for
$$\dd_k:=  \Big(\ff 1 2+ \ff {4^{q-1}(ck)^q \e^{2^{q-1}k^q t} }{q\ll+ 2^{q-1}k^q}\Big)<1, \ \ k\in (0,k_p).$$
So,   \eqref{Q'} holds for some constant $c'>0$ due to the semigroup property $P_{t+s}^*= P_t^*P_s^*.$
\end{proof}

To verify   condition \eqref{LNA}, we present below a Harris type theorem  on the uniform ergodicity for a family of Markov processes.

\beg{lem}\label{L4} Let $(E,\rr)$ be a metric space and  let $\{(P_t^i)_{t\ge 0}: i\in I\}$ be a family of Markov semigroups on $\B_b(E)$. If there exist $t_0,t_1>0$ and measurable set $B\subset E$ such that
\beq\label{PST}\aa:= \inf_{i\in I, x\in E} P_{t_0}^i 1_B(x)>0,\end{equation}
\beq\label{PST2} \bb:=\sup_{i\in I, x,y \in B} \|(P_{t_1}^i)^*\dd_x- (P_{t_1}^i)^*\dd_y\|_{var} <2,\end{equation}
then there exists $c>0$ such that
\beq\label{UEE} \sup_{i\in I, x,y \in E} \|(P_{t}^i)^*\dd_x- (P_{t}^i)^*\dd_y\|_{var} \le c \e^{-\ll t},\ \ t\ge 0\end{equation}
holds for $\ll:=\ff 1{t_0+t_1} \log\ff 2{2-\aa^2(2-\bb)}>0.$
\end{lem}

\beg{proof} The proof is more or less standard. By the semigroup property, we have
\beg{align*} &\|(P_{t_0+t_1}^i)^*\dd_x- (P_{t_0+t_1}^i)^*\dd_y\|_{var}\\
&=\sup_{|f|\le 1} \bigg|\int_{E\times E}  \big(P_{t_1}^if(x')- P_{t_1}^if(y')\big) \{(P_{t_0}^i)^*\dd_x\}(\d x') \{(P_{t_0}^i)^*\dd_y\}(\d y')\bigg| \\
&\le \int_{B\times B} \|(P_{t_1}^i)^*\dd_{x'}- (P_{t_1}^i)^*\dd_{y'}\|_{var}\{(P_{t_0}^i)^*\dd_x\}(\d x') \{(P_{t_0}^i)^*\dd_y\}(\d y') \\
&\quad + 2 \int_{(B\times B)^c}  \{(P_{t_0}^i)^*\dd_x\}(\d x') \{(P_{t_0}^i)^*\dd_y\}(\d y') \\
&\le \bb \{P_{t_0}^i 1_B(x)\} P_{t_0}^i 1_B(y)+ 2 \big[1-\{P_{t_0}^i 1_B(x)\} P_{t_0}^i 1_B(y)\big]\le 2-\aa^2(2-\bb).\end{align*}
Thus, for $\dd:= \ff{2-\aa^2(2-\bb)} 2 <1$, we have
$$\|(P_{t_0+t_1}^i)^*\dd_x- (P_{t_0+t_1}^i)^*\dd_y\|_{var}\le \dd \|\dd_x-\dd_y\|_{var},\ \ x,y\in E.$$
Combining this with the semigroup property,  we find constants $c>0$ such that \eqref{UEE} holds for the claimed $\ll>0$.
\end{proof}

\beg{proof}[ Proof  Theorem  \ref{T3}]
 According to Theorems \ref{T1.6.1}-\ref{T2} and Lemma \ref{L3},   it suffices to prove  \eqref{LNA} for small $k>0$. By Lemma \ref{L4}, we only need to prove
\eqref{PST} and \eqref{PST2} for the family $\{P_t^\gg: \gg\in \scr P\},$ where $P_t^\gg f(x):= \int_{\bar D} f\d(P_t^\gg)^*\dd_x.$

(a)  Proof of \eqref{PST2}. Let us fix $\gg\in \scr P$, and let $X_t^{x,\gg}$ solve \eqref{EGG} with $X_0^\gg=x$. For any $\nu\in\scr P$, by Girsanov's theorem we have
$$P_t^{\nu}f(x)= \E[f(X_t^{x,\gg}) R_t^{x,\gg,\nu}],\ \ t\ge 0,$$
where
$$R_t^{x,\gg,\nu}:=\e^{\int_0^t\<\zeta(\gg,\nu)(X_s^{x,\gg}),\d W_s\> -\ff 1 2|\zeta(\gg,\nu)(X_s^{x,\gg})|^2\d s}.$$
By {\bf (H2)}, Girsanov's theorem  and Pinsker's inequality,  we obtain 
$$\|(P_t^{\gg})^* \dd_z-(P_t^{\nu})^*\dd_z\|_{var}^2 \le  (\E|R_t^{\gg,\nu}-1|)^2\le k^2 t^{\ff 2 q}  \|\gg-\nu\|_{var}^2\le 4 k^2t^{\ff 2 q},\ \ t\ge 0, z\in\bar D, \nu\in\scr P.$$
Taking $t=t_1= (4k)^{-q},$  we obtain
\beq\label{PQ1} \sup_{\nu\in \scr P}\|(P_{t_1}^{\gg})^* \dd_z-(P_{t_1}^{\nu})^*\dd_z\|_{var}\le \ff 1 2,\ \ z\in\bar D,\nu\in \scr P.\end{equation}
On the other hand, by \eqref{STF},  there exists $x_0\in D$ and a constant $\vv>0$ such that  $B(x_0,\vv)\subset D$ and
$$\|(P_{t_1}^\gg)^*\dd_x- (P_{t_1}^\gg)^*\dd_y\|_{var} \le \ff 1 4,\ \ x,y\in B(x_0,\vv).$$
Combining this with \eqref{PQ1} we derive
$$\sup_{\nu\in \scr P} \|(P_{t_1}^\nu)^*\dd_x- (P_{t_1}^\nu)^*\dd_y\|_{var} \le \ff 3 2<2,\ \ x,y\in B(x_0,\vv).$$
So,  \eqref{PST2} holds for $B=B(x_0,\vv)$.

(b)   Let $u$ solve \eqref{*1}  for $f=-b^{(0)}$ and large $\ll>0$ such that  \eqref{*4} holds, and let $\Theta(x)=x+ u(x)$. By {\bf (H1)}, 
we see that     \eqref{GMM}  holds with  $Y_t^{x,\nu}:= \Theta(X_t^{x,\nu})$   replacing $Y_t$ for  
all $\nu\in \scr P$.  So, by $H(\infty)<\infty$ and the argument leading to \eqref{GMM2}, we obtain
$$\sup_{\nu\in \scr P, x\in\bar D} \E  [V(Y_t^{x,\nu})] \le \theta^{-1} k,\ \ t\ge kH(\infty)=:t_2.$$
This together with \eqref{*GMM} implies
$$\sup_{\nu\in \scr P, x\in \bar D} \E [V(X_t^{x,\nu})]  \le \theta^{-2} k,\ \ t\ge t_2.$$
Letting ${\bf K}:=\{V\le 2\theta^{-2} k\}$, we derive
\beq\label{ESTN} \inf_{\nu\in \scr P, x\in\bar D} P_{t_2}^\nu 1_{\bf K}(x) \ge \ff 1 2.\end{equation}

On the other hand,  by  Girsanov's theorem and Schwartz's inequality, we find a constant $c_0>0$ such that
$$P_{1}^\nu 1_{B(x_0,\vv)}(x) = \E \big[1_{B(x_0,\vv)}(X_1^{x,\gg}) R_{1}^{x,\gg,\nu}\big]
\ge \ff{\{\E  1_{B(x_0,\vv)}(X_{1}^{x,\gg}) \}^2}{\E R_{1}^{x,\gg,\nu}}\ge c_0 (P_{1}^\gg 1_{B(x_0,\vv)}(x))^2.$$
Since $\bf K$ is bounded, combining this with Lemma \ref{L2} for $P_t^\gg$, we find a constant $c_1>0$ such that
$$\inf_{\nu\in \scr P, x\in {\bf K}} P_{1}^\nu 1_{B(x_0,\vv)}(x) \ge c_1.$$ This together with    \eqref{ESTN} and the semigroup  property  yields
$$P_{t_2+1}^\nu 1_{B(x_0,\vv)}(x) \ge P_{t_2 }^\nu \big\{1_{\bf K} P_1^\nu 1_{B(x_0,\vv)}\}(x)\ge c_1 P_{t_2 }^\nu1_{\bf K}(x) \ge \ff{c_1}2>0,\ \ x\in\bar D,\nu\in \scr P.$$
Therefore, \eqref{PST} holds for $t_0= t_2+1.$ \end{proof}

\section{Exponential ergodicity  in   entropy and   Wasserstein distance}

  In this section, we consider the following reflecting McKean-Vlasov SDE where the noise may also be distribution dependent:
\beq\label{E1N} \d X_t= b_t(X_t,\L_{X_t}) \d t+ \si_t(X_t,\L_{X_t})\d W_t + \n(X_t) \d l_t,\ \ t\ge 0,\end{equation}
where 
$$b: [0,\infty)\times \bar D\times\scr P\to\R^d,\ \ \si: [0,\infty)\times \bar D\times\scr P\to\R^d\otimes\R^n$$ are measurable. 
We study the  exponential ergodicity  under entropy and weighted Wasserstein distance for  dissipative and partially dissipative cases respectively, such that 
the corresponding results in \cite{20RWb} and \cite{W21a} are extended to the reflecting setting.  
 For simplicity, we  only consider convex $D$,  for which the local time on boundary does not make trouble in the study.

\emph{\beg{enumerate} \item[{\bf (A7)}] Let $k> 1$, $\scr P_k :=\{\mu\in \scr P: \mu(|\cdot|^k)<\infty\}.$   $D$ is convex, 
    $b$ and $ \si$ are bounded on bounded subsets of  $[0,\infty)\times\bar D\times\scr P_k(\bar D)$, and  the following two conditions hold.
\item[$(1)$] For any $T>0$ there exists a constant  $K>0$ such that
  \beg{align*} &  \|\si_t(x,\mu)-\si_t(y,\nu)\|^2_{HS}+2\<x-y,b_t(x,\mu)-b_t(y,\nu)\>^+ \\
  &\le K  \big\{|x-y|^2+|x-y|\W_k(\mu,\nu)+1_{\{k\ge 2\}}\W_k(\mu,\nu)^2\big\},\ \ t\in [0,T], x,y\in \bar D, \mu,\nu\in \scr P_k(\bar D).\end{align*}
\item[$(2)$] There exists a subset $\tt\pp D\subset \pp D$ such that
\beq\label{A21} \<y-x, \n(x)\>\ge 0,\ \ x\in \pp D\setminus \tt\pp D,\  y\in \bar D,\end{equation} and when $\tt\pp D\ne \emptyset,$ there exists $\tt\rr\in C_b^2(\bar D)$ such that $\tt\rr|_{\pp D}=0$,
$\<\nn\tt\rr, \n\>|_{\pp D}\ge 1_{\tt\pp D}$ and
\beq\label{A22}  \sup_{(t,x)\in   [0,T]\times \bar D} \big\{\|(\si_t^\mu)^*\nn\tt\rr\|^2(x)+ \<b_t^\mu,\nn \tt\rr\>^-(x)\big\}<\infty,\ \ \mu\in C([0,T]; \scr P_k(\bar D)).\end{equation} \end{enumerate}}
According to \cite{W21b}, this assumption implies the well-posedness of \eqref{E1N} for distributions in $\scr P_k$. Let
$P_t^*\mu=\L_{X_t}$ for the solution with $\L_{X_0}=\mu\in \scr P_k$.

\subsection{Dissipative case: exponential convergence in entropy and $\W_2$}

In this part, we study  the  exponential ergodicity of $P_t^*$ in entropy and $\W_2$. For probability measures $\mu_1,\mu_2$ on $\bar D$, let
$$\Ent(\mu_1|\mu_2):=\beg{cases} \mu_2(f\log f),\ &\text{if}\ \d\mu_1=f\d\mu_2,\\
\infty,\ &\text{if}\ \mu_1\ \text{is\ not\ absolutely\ continuous\ w.r.t.\ }\mu_2\end{cases}$$
be the relative entropy of $\mu_1$ w.r.t. $\mu_2$, and let
$$\W_2(\mu_1,\mu_2):=\inf_{\pi\in \C(\mu_1,\mu_2)} \bigg(\inf\int_{\bar D\times\bar D} |x-y|^2\pi(\d x,\d y)\bigg)^{\ff 1 2} 
$$ be the quadratic Wasserstein distance, where $\C(\mu_1,\mu_2)$ is the set  of  all couplings for $\mu_1$ and $\mu_2$.  The following result extends  the corresponding one derived   in   \cite{20RWb} for  McKean-Vlasov SDEs without reflection. 

\beg{thm}\label{T6} Let $D$ be convex and   $(\si,b) $ satisfy {\bf (A7)} with $k=2$.
Let $K_1, K_2\in  L_{loc}^1([0,\infty);\R) $     such that
\beq\label{DSS} \beg{split} &2\<b_t(x,\mu)-b_t(y,\nu), x-y\> + \|\si_t(x,\mu)-\si_t(y,\nu)\|_{HS}^2\\
&\le K_1(t) |x-y|^2+ K_2 (t) \W_2(\mu,\nu)^2,\ \ t\ge 0.\end{split} \end{equation}
Then  and $P_t^*$ satisfies
\beq\label{NNW} \W_2(P_{t}^*\mu, P_{t}^*\nu)^2 \le \e^{\int_0^t (K_1+K_2)(r)\d r}\W_2(\mu,\nu)^2,\ \  \mu,\nu\in \scr P_2(\bar D), \ t\ge  0. \end{equation}  Consequently, if $(b_t,\si_t)$ does not depend on $t$ and $\ll:=-(K_1+K_2)>0$,
then $P_{t}^*$ has a unique invariant probability measure $\bar\mu$  such that
\beq\label{BBE} \W_2(P_{t}^*\mu, \bar\mu)^2 \le \e^{-\ll t}\W_2(\mu,\bar\mu)^2,\ \  \mu\in \scr P_2(\bar D), \ t\ge 0,\end{equation}
and the following assertions hold: 
\beg{enumerate} \item[$(1)$] When  $\si_t(x,\mu)=\si_t(x)$ does not depend on $\mu$ and $\si\si^*$ is invertible with $\|\si\|_\infty+\|(\si\si^*)^{-1}\|_\infty<\infty$,   there exists a constant $c>0$ such that
\beq\label{ET*} \Ent(P_t^*\mu|\bar\mu) \le c\e^{-\ll t} \W_2(\mu,\bar\mu)^2,\ \ t\ge 1, \mu\in \scr P_2(\bar D).\end{equation}
\item[$(2)$] When   $\si(x,\mu)=\si(\mu)$ does not depend on $x$,   there exists a constant $c>0$ such that $\bar\mu$ satisfies the following log-Sobolev inequality and  Talagrand inequality:
    \beq\label{LSI}  \bar\mu(f^2\log f^2)\le c\bar\mu(|\nn f|^2),\ \ f\in C_b^1(\R^d), \bar\mu(f^2)=1,\end{equation}
\beq\label{TTI} \W_2(\mu,\bar\mu)^2\le c \Ent(\mu|\bar\mu),\ \ \mu\in \scr P_2.\end{equation}
\item[$(3)$] When   $\si(x,\mu)=  \si$   is constant with $\si\si^*$ invertible,   there exists a constant $c>0$ such that
\beq\label{NNW2} \beg{split} &\W_2(P_t^*\mu, \bar\mu)^2+\Ent(P_t^*\mu|\bar \mu)\\
&\le c \e^{-\ll t} \min\big\{\W_2(\mu,\bar\mu)^2,\Ent(\mu|\bar \mu)\big\},\ \ t\ge 1,\mu\in\scr P_2(\bar D).\end{split} \end{equation}\end{enumerate}
\end{thm}

\beg{proof}   The well-posedness is ensured by \cite[Theorem  3.3]{W21a}.  Since $D$ is convex, by Remark 2.1 in \cite{W21a},
\beq\label{CVX} \<y-x, \n(x)\>\ge 0,\ \ y\in\bar D, x\in\pp D, \n(x)\in \scr N_x.\end{equation}  For any $\mu,\nu\in \scr P_2(\bar D)$, let
$X_0^\mu$ and $X_0^\nu$ be $\F_0$-measurable such that
\beq\label{NAB} \L_{X_{0}^\mu}=\mu, \ \ \L_{X_{0}^\nu}=\nu,\ \ \E|X_{0}^\mu-X_{0}^\nu|^2=\W_2(\mu,\nu)^2.\end{equation}
  By \eqref{DSS}, \eqref{CVX},  and   applying It\^o's formula to $|X_{t}^\mu-X_{t}^\nu|^2,$ where $(X_{t}^\mu)_{t\ge 0}$ and $(X_{t}^\nu)_{t\ge 0}$ solve \eqref{E1N}, we  obtain
$$\d |X_t^\mu-X_t^\nu|^2\le \big\{ K_1(t)|X_t^\mu-X_t^\nu|^2 +K_2(t)\W_2(\P_t^*\mu, P_t^*\nu)^2\big\}\d t +\d M_t$$
for some martingale $M_t$. Combining this with \eqref{NAB}, $\W_2(P_t^*\mu,P_t^*\nu)^2\le \E|X_t^\mu-X_t^\nu|^2$, and Gronwall's lemma, we prove
 \eqref{NNW}.

 Let $(b_t,\si_t)$ do not depend on $t$ and $\ll:=-(K_1+K_2)>0$. Then \eqref{NNW} implies the uniqueness of $P_t^*$-invariant probability measure $\bar\mu\in \scr P_2(\bar D)$ and
\eqref{BBE}.

The existence of $\bar\mu$   follows from a standard argument by showing that for $x_0\in D$, $\{P_t^*\dd_{x_0}\}_{t\ge 0}$ is a $\W_2$-Cauchy family as $t\to\infty$. Since the term of local time does not make trouble due to \eqref{CVX},
the proof is completely similar to that of \cite[Theorem 3.1]{W18} for the case $D=\R^d$, so we skip the details to save space. Below we prove statements (1)-(3) respectively. 

(1) When $\si_t(x,\mu)=\si_t(x)$  and $\si\si^*$ is invertible with $\|\si\|_\infty+\|(\si\si^*)^{-1}\|_\infty<\infty$,  by Theorem 4.2 in \cite{W21a}, {\bf (A7)} with $k=2$ implies the log-Harnack inequality
$$\Ent(P_1^*\mu|\bar\mu)\le c \W_2(\mu,\bar\mu)^2,\ \ \mu\in \scr P_2(\bar D)$$
for some constant $c>0$. So, \eqref{ET*} follows from \eqref{BBE} and   $P^*_t= P_1^*P_{t-1}^*$ for $t\ge 1$.

 (2) Let  $\si(x,\mu)=  \si(\mu)$ be independent of  $x$. Consider the SDE
\beq\label{SS4} \d\bar X_{t}^x= b(\bar X_{t}^x, \bar\mu)\d t +  \si(\bar\mu) \d W_t+\n(\bar X_t^x)\d l_t,\ \ t\ge s, \bar X_{0}^x=x\in\bar D.\end{equation}
The associated Markov semigroup $\{\bar P_{t}\}_{t\ge   0}$  is  given by
$$\bar P_{t}  f(x):= \E f(\bar X_{t}^x),\ \ t\ge   0, f\in \B_b(\bar D), x\in\bar D.$$  Let  $\bar P_{t}^*$ be given by
$$(\bar P_{t}^*\mu)(f):= \mu(\bar P_{t}f),\ \ \mu\in \scr P, t\ge   0, f\in \B_b(\bar D).$$
 Since \eqref{DSS} with $x=y$ implies $K_2\ge 0$, we have
\beq\label{KK1} K_1  \le -\ll<0.\end{equation}
As explained in the above proofs of \eqref{NNW} and \eqref{BBE},  this implies that  $\bar P_t^*$ has a unique invariant probability measure
$ \tt\mu$  such that
\beq\label{WCC} \lim_{t\to\infty} \bar P_t f(x)=\tt\mu(f),\ \ f\in C_b(\bar D), x\in\bar D.\end{equation}
Since $\bar\mu$ is the unique invariant probability measure of $P_t^*$,    and  when the initial distribution is $\bar\mu,$ the SDE \eqref{SS4} coincides with \eqref{E1N},   we conclude that $\tt\mu=\bar \mu$.
Hence, \eqref{WCC} yields
\beq\label{RTT} \bar\mu(f)=\lim_{t\to\infty} P_{t} f(x_0),\ \ f\in C_b(\bar D), x_0\in D.\end{equation}
Now, by It\^o's formula, \eqref{CVX}  and \eqref{DSS} with $(b_t,\si_t)$ independent of $t$, we   obtain
$$|\bar X_{t}^x-\bar X_{t}^y|^2\le \e^{K_1t} |x-y|^2,\ \ x, y\in \bar D, t\ge   0.$$
This and  \eqref{KK1}   imply
\beq\label{GFF} \beg{split} &|\nn \bar P_{t}f(x)|:=\limsup_{y\to x} \ff{|\bar P_{t}f(x)-\bar P_{t}f(y)|}{|x-y|}
\le \limsup_{y\to x} \ff{\E|f(\bar X_t^x)-f(\bar X_t^y)|}{|x-y|} \\
&\le \e^{-\ff{\ll t}2} \limsup_{y\to x} \E\ff{|f(\bar X_t^x)-f(\bar X_t^y)|}{|\bar X_t^x-\bar X_t^y|}= \e^{-\ll t/2} \bar P_t |\nn f|(x),\ \ t\ge 0, f\in C_b^1(\bar D).\end{split}\end{equation}
On the other hand, we have
$$ \pp_t\bar P_{t}f=   \bar L \bar P_{t}f,\ \    \<\n, \nn \bar P_{t} f\>|_{\pp D}=0,\ \ t\ge 0, f\in C_N^2(\bar D),$$
where $C_N^2(\bar D)$ is the set of $f\in C_b^2(\bar D)$ satisfying with $\<\n,\nn f\>|_{\pp D}=0,$ and
$$\bar L:= \ff 1 2 {\rm tr} \{(\bar \si\bar \si^*)  \nn^2\}+\nn_{b(\cdot,  \bar\mu)},\ \ \bar\si:=\si(\bar\mu),\ \ s\ge 0.$$
So, by It\^o's formula, for any $\vv>0$ and $f\in C_N^2(\bar D)$,
$$ \d \big\{( \bar P_{t-s} (\vv+f^2)) \log  \bar P_{t-s} (\vv+f^2)\big\}(\bar X_s) = \Big\{\ff{|\bar\si^*\nn \bar P_{t-s} f^2|^2}{\vv+ \bar P_{t-s} f^2} \Big\}\d t +\d M_s^\vv,\ \ s\in [0,t]$$
holds for some martingale $(M_s^\vv)_{s\in [0,t]}$.  Combining this with \eqref{GFF}, we    find a constant $c>0$ such that for any $f\in C_N^2(\R^d)$,
\beg{align*} &\bar P_t\big\{ (\vv+f^2 ) \log  (\vv+ f^2)\big\} -  (\vv+\bar P_{t} f^2)\log (\vv+ \bar P_{t} f^2 )  \\
&= \int_0^t \bar P_s \ff{|\bar \si^* \nn \bar P_{t-s} f^2|^2}{\vv+ \bar P_{t-s}f^2}  \d s \le 4 (c_1 \|\bar\si\|_\infty)^2 \int_0^t \e^{-\ll (t-s)} \bar P_s\bar P_{t-s} |\nn f|^2 \d s\\
&= 4 (c_1 \|\bar\si\|_\infty)^2(\bar P_t|\nn f|^2)  \int_0^t \e^{-\ll (t-s)}   \d s \le c \bar P_t|\nn f|^2,\ \ t\ge 0,\vv>0.\end{align*}
By   letting first $\vv\downarrow 0$ then $ t\to\infty$, we deduce from this and  \eqref{RTT} that
$$ \bar\mu(f^2\log f^2)\le c_2 \bar\mu(|\nn f |^2),\ \ f\in C_N^2(\bar D), \bar\mu(f^2)=1$$ holds for some constant $c_2>0$.
This implies \eqref{LSI} by an approximation argument, indeed the inequality holds for $f\in H^{1,2}(\bar\mu)$ with $\bar\mu(f^2)=1.$
 According to Lemma \ref{LTI} below,  \eqref{TTI} holds.

(3) Let $\si$ be constant with $\si\si^*$ invertible. Then    \eqref{NNW2} follows from   \eqref{BBE}, \eqref{ET*}  and \eqref{TTI}. \end{proof}

The following result on the Talagrand inequality  is known by   \cite{BGL} when $\bar\mu(\d x)= \e^{V(x)}\d x $ for some $V\in C(\R^d),$
 which is first proved in \cite{OV} on Riemannian manifolds under a curvature condition, see also \cite{W04} for more general results.
 We extend it to general probability measures for the above application to $\bar \mu$  which is  supported  on $\bar D$ rather than $\R^d$.

\beg{lem}\label{LTI}  Let $c>0$ be a constant and $\bar \mu\in \scr P_2(\R^d).$  Then the log-Sobolev inequality
\eqref{LSI}  implies \eqref{TTI}.
\end{lem}

\beg{proof} By an approximation argument, we only need to prove for $\mu=\varrho \bar\mu$ for some  density $\varrho \in C_b(\R^d)$
Let $P_t^{(0)}$ be the Ornstein-Uhlenbeck semigroup generated by
$\DD-x\cdot\nn$ on $\R^d.$ We have
$$|\nn P_t^{(0)} f|\le P_t^{(0)}|\nn f|,\ \ P_t^{(0)}(f^2\log f^2)\le t P_{t}^{(0)}|\nn f |^2 +(P_t^{(0)}f^2)\log P_t^{(0)}f^2,\ \ f\in C_b^1(\R^d).$$
Combining this with \eqref{LSI}, we see that $\bar \mu_t:= (P_t^{(0)})^*\bar\mu$ satisfies
\beg{align*} &\bar\mu_t (f^2\log f^2) = \bar\mu(P_t^{(0)}(f^2\log f^2)) \le t \bar\mu_t(|\nn f|^2) + \bar\mu((P_t^{(0)} f^2)\log P_t^{(0)}f^2)\\
&\le t \bar\mu_t(|\nn f|^2) + c\bar\mu\Big(\Big|\nn \ss{P_t^{(0)} f^2}\Big|^2\Big) + \bar\mu_t(f^2)\log \bar\mu_t(f^2)\\
&\le (t+ c)\bar\mu_t(|\nn f|^2) +  \bar\mu_t(f^2)\log \bar\mu_t(f^2),\ \ f\in C_b^1(\R^d),\ \ t>0,\end{align*}
where the last step follows from the gradient estimate $|\nn P_t^{(0)}f|\le P_t^{(0)}|\nn f|,$ which and  the Schwarz inequality imply
$$\Big|\nn \ss{P_t^{(0)} f^2}\Big|^2=\ff {|\nn P_t^{(0)} f^2|^2}{ 4 P_t^{(0)}f^2}\le \ff{\{P_t^{(0)}(|f\nn f|)\}^2}{P_t^{(0)}f^2}\le P_t^{(0)}|\nn f|^2.$$
Therefore, $\bar\mu_t$ satisfies the log-Sobolev inequality with constant $t+c$ and has smooth strictly positive density. According to \cite{BGL}, we have
$$\W_2(\mu,\bar\mu_t)^2\le (t+c) \Ent(\mu|\bar\mu_t),\ \ \mu\in \scr P_2(\R^d).$$
Since $\W_2(\bar\mu_t, \bar\mu)\to 0$ as $t\to 0$,   and  $\mu=\varrho \bar\mu$ with $  \varrho\in C_b(\R^d)$, this implies
\beg{align*} &\W_2(\mu,\bar\mu)^2= \lim_{t\downarrow 0} \W_2(\mu,\bar\mu_t)^2\le  \lim_{t\downarrow 0} (t+c) \Ent(\mu|\bar\mu_t)\\
&=  \lim_{t\downarrow 0} (t+c) \bar\mu((P_t^{(0)} \varrho)\log P_t^{(0)} \varrho)
 =c\bar\mu( \varrho\log \varrho).\end{align*}
Therefore,   \eqref{TTI} holds.
\end{proof}

\subsection{Partially dissipative case: exponential convergence in $\W_\psi$}

In this part, we consider the partially dissipative case such that \cite[Theorem 3.1]{W21a} is extended to the reflecting setting.
For  any $\kk>0$, let
\beq\label{PKK} \beg{split}
  \Psi_\kk:= \big\{\psi\in C^2((0,\infty))\cap C^1([0,\infty)):\ & \psi(0)=0, \ \psi'|_{(0,\infty)}>0, \|\psi'\|_\infty<\infty \\
  &\ r\psi'(r)+ r^2\{\psi''\}^+(r) \le \kk \psi(r)\ \text{for}\ r> 0 \big\}.\end{split}\end{equation}

For $\psi\in \Psi_\kk$, we introduce the associated  Wasserstein $``$distance'' (also called transportation cost)
 \beq\label{PKK'} \W_\psi (\mu,\nu):=   \inf_{\pi\in\C(\mu,\nu)} \int_{\bar D\times\bar D} \psi(|x-y|) \pi(\d x,\d y),\ \ \mu,\nu\in \scr P_\psi.\end{equation}
Then $\W_\psi$ is a complete quasi-metric on the space
$$\scr P_{\psi}:=\big\{\mu\in \scr P: \mu(\psi(|\cdot|))<\infty\big\}.$$

\emph{\beg{enumerate} \item[{\bf (A8)}]     $\si_t(x,\mu)=\si_t(x)$ does not depend on $\mu$ so that \ref{E1N} reduces to \eqref{E1}.
\item[$(1)$] $($Ellipticity$)$  There exist    $\aa\in C([0,\infty); (0,\infty))$ and
$\hat \si\in \B( [0,\infty)\times \bar D; R^d\otimes\R^d)$ such that
$$\si_t(x)\si_t(x)^*=\aa_t  {\bf I}_d+ \hat \si_t(x)\hat \si_t(x)^*,\ \ t\ge 0, x\in\bar D.$$
\item[$(2)$] $($Partial dissipativity$)$ Let  $\psi\in \Psi_\kk$ in \eqref{PKK} for some $\kk>0$, $\gg\in C([0,\infty))$ with $\gg(r)\le K r$ for some constant $K>0$ and all $r\ge 0$,   such that
\beq\label{A2E}  2\aa_t \psi''(r) +(\gg\psi')(r)\le -\zeta_t \psi(r),\ \ r\ge 0,t\ge 0\end{equation} holds for some for some $\zeta\in C([0,\infty);\R).$
Moreover, $b\in C([0,\infty)\times \bar D\times \scr P_\psi),$   and   there exists     $\theta \in C([0,\infty); [0,\infty))$       such that
\beq\label{A3E} \beg{split}  & \<b_t(x,\mu)-b_t(y,\nu), x-y\> +\ff 1 2 \|\hat \si_t(x)-\hat \si_t(y)\|_{HS}^2\\
&\quad \le |x-y|  \big\{\theta_t  \W_\psi (\mu,\nu)  + \gg(|x-y|)\big\},\ \ t\ge 0, x,y\in \bar D, \mu,\nu\in \scr P_\psi.   \end{split}\end{equation}\end{enumerate}}

 \beg{thm} \label{T7} Let $D$ be convex and assume  {\bf (A8)}, where   $\psi''\le 0$ if $\hat \si$ is non-constant. Then   $\eqref{E1}$ is well-posed for distributions in
 $\scr P_\psi$, and   $P_t^*$ satisfies
 \beq\label{EXP1'}  \W_\psi(P_t^*\mu, P_t^*\nu)\le  \e^{-\int_0^t\{\zeta_s-\theta_s\|\psi'\|_\infty\}\d s}  \W_\psi(\mu,\nu),\ \ t\ge 0, \mu,\nu \in \scr P_\psi.\end{equation}
 Consequently,
 if    $(b_t,\si_t,\zeta_t,\theta_t)$ do  not depend on $t$ and $\zeta>\theta \|\psi'\|_\infty$,    then $P_t^*$ has a unique invariant probability measure
$\bar\mu\in \scr P_\psi$ such that
\beq\label{EXP2'}  \W_\psi(P_t^*\mu, \bar\mu)\le  \e^{-(\zeta-\theta \|\psi'\|_\infty) t}  \W_\psi(\mu,\bar\mu),\ \ t\ge 0, \mu \in \scr P_\psi.\end{equation}
  \end{thm}

\beg{proof}  Since $D$ is convex, the proof is similar to that of \cite[Theorem 3.1]{W21a}. We outline it below for complement.

By Theorem \ref{T3},  the well-posedness follows from {\bf (A8)}(1) and {\bf (A8)}(2).
Next, according to the proof of Theorem \ref{T6}(2) with $\W_\psi$ replacing $\W_2$, the second assertion follows from the first. So, in the following we only prove \eqref{EXP1'}.

  For any $s\ge 0$, let $(X_s,Y_s)$ be $\F_s$-measurable such that
 \beq\label{O1}\L_{X_s}=P_s^*\mu, \ \ \L_{Y_s}=P_s^*\nu,\ \  \W_\psi(P_s^*\mu,P_s^*\nu)= \E\psi(|X_s-Y_s|).\end{equation}
 Let $W_t^{(1)}$ and $ W_t^{(2)}$ be two independent $d$-dimensional Brownian motions and consider the following
 SDE:
 \beq\label{E*A} \d X_t= b_t(X_t, P_t^*\mu)\d t  + \ss{\aa_t}\d W_t^{(1)}+ \hat \si_t(X_t)\d W_t^{(2)}+\n(X_t)\d l_t^X,\ \ t\ge s,\end{equation}
 where $l_t^X$ is the local time of $X_t$ on $\pp D$.
  By  Theorem \ref{T6}, {\bf (A8)}(1) and {\bf (A8)}(2) imply that  this SDE. By  $\si_t\si_t^* = \aa_t  {\bf I}_d+ \hat \si_t \hat \si_t^*$, we have
  \beg{align*} &\si_t^*(\si_t\si_t^*)^{-1}\big\{\aa_t+\hat\si_t  \hat\si_t^*\big\} (\si_t\si_t^*)^{-1}\si_t+ \big\{{\bf I}_m - \si_t^*(\si_t\si_t^*)^{-1}\si_t\big\}^2\\
  &= \si_t^*(\si_t\si_t^*)^{-1} \si_t+{\bf I}_m-  \si_t^*(\si_t\si_t^*)^{-1} \si_t={\bf I}_m.\end{align*} So, for an $m$-dimensional Brownian motion $W^{(3)}$ independent of $(W^{(1)}, W^{(2)}$,
  $$W_t:= \int_0^t\big\{\si_s^*(\si_s\si_s^*)^{-1}\big\}(X_s)  \big\{\ss{\aa_s}\d W_s^{(1)}+ \hat \si_s(X_s)\d W_s^{(2)}\big\}+\int_0^t \\big\{I_m - \si_s^*(\si_s\si_s^*)^{-1}\si_s\big\}(X_s)\d W_s^3$$
  is an $m$-dimensional Brownian motion such that
  $$\si_t(X_t)\d W_t=  \ss{\aa_t}\d W_t^{(1)}+ \hat \si_t(X_t)\d W_t^{(2)}.$$
  Thus, by   the weak uniqueness of \eqref{E1}, we have    $\L_{X_t}= P_{s,t}^* P_s^*\mu=P_t^*\mu, t\ge s,$ where for $\gg\in \scr P_\psi$ we denote $P_{s,t}^*\gg=\L_{X_t}$ for
  $X_t$ solving \eqref{E*A} with $\L_{X_s}=\gg.$

To construct the   coupling with reflection, let
   $$u(x,y)=   \ff{x-y}{|x-y|},\ \ x\ne y\in \R^d.$$
 We consider the SDE for $t\ge s$:
 \beq\label{E*1}  \d Y_t= b_t(Y_t, P_t^*\nu)\d t  + \ss{\aa_t}\big\{ {\bf I}_d-2u(X_t,Y_t)\otimes u(X_t,Y_t) 1_{\{t<\tau\}}\big\} \d W_t^{(1)}+ \hat \si_t(Y_t)\d W_t^{(2)}  +\d l_t^Y,\end{equation}
  where
  $$\tau:= \inf\{t\ge s: Y_t=X_t\}$$ is the coupling time.
 Since the coefficients in noises are  Lipschitz continuous outside a neighborhood of the diagonal,  by  \cite[Theorem 1.1]{Hino},
 \eqref{E*1} has a unique solution up to the coupling time $\tau$.
 When $t\ge \tau$, the equation of $Y_t$ becomes
 \beq\label{E*1'} \d Y_t= b_t(Y_t, P_t^*\nu)\d t  + \ss{\aa_t}\d W_t^{(1)}+ \hat \si_t(Y_t)\d W_t^{(2)}+\d l_t^Y,\end{equation} which is well-posed under  {\bf (A8)}(1) and {\bf (A8)}(2) according to Theorem \ref{T3}.
 So,  \eqref{E*1} is well-posed and    $\L_{Y_t}= P_t^*\nu $  by the same reason leading to $\L_{X_t}= P_t^*\mu.$
 Since $D$ is convex, \eqref{CVX} holds. So, by  {\bf (A8)}(1) and {\bf (A8)}(2) for $\psi\in \Psi$ with $\psi''\le 0$ when $\hat\si_t$ is non-constant, and applying It\^o's formula, we obtain
\beq\label{ITP'}  \beg{split} \d \psi(|X_t-Y_t|)  \le & \big\{\theta_t \psi'(|X_t-Y_t|)   \W_\psi(P_t^*\mu,P_t^*\nu) -\zeta_t  \psi(|X_t-Y_t|)
   \big\}\d t   \\
& +  \psi'(|X_t-Y_t|)\Big[2 \ss{\aa_t}   \Big\<u(X_t,Y_t), \d W_t^{(1)}\Big\> \\
&+  \Big\<u(X_t,Y_t), (\hat \si_t(X_t)-\hat \si_t(Y_t))\d W_t^{(2)}\Big\>\Big],\ \  s\le t<\tau.\end{split}\end{equation}
By a standard argument and noting that $\psi(|X_{t\land\tau}, Y_{t\land\tau}|) 1_{\{\tau\le t\}}=0$,  this implies
\beg{align*} & \e^{\int_s^t \zeta_p\d p}   \E\big[\psi(|X_{t\land\tau}-Y_{t\land\tau}|)\big]= \E\big[\e^{\int_s^{t\land\tau} \zeta_p\d p} \psi(|X_{t\land\tau}-Y_{t\land\tau}|)\big]\\
&\le \E \psi(|X_s-Y_s|) + \|\psi'\|_\infty \int_s^{t\land\tau} \theta_r\e^{\int_s^{r}\zeta_p\d p }\W_\psi(P_r^*\mu, P_r^*\nu)\d r,\ \ t\ge s.\end{align*}
Consequently,
\beq\label{W1'} \beg{split} &\E \psi(|X_{t\land\tau}-Y_{t\land\tau}|)\\
&\le \e^{-\int_s^t\zeta_r\d r}\E \psi(|X_s-Y_s|) + \|\psi'\|_\infty \int_s^{t\land\tau} \theta_r\e^{-\int_r^{t}\zeta_p\d p }\W_\psi(P_r^*\mu, P_r^*\nu)\d r,\ \ t\ge s.\end{split}\end{equation}
On the other hand, when $t\ge \tau$, by {\bf (A8)}(2) and applying
It\^o's formula for \eqref{E*A} and \eqref{E*1'},   we find a constant $C>0$ such that
\beg{align*} \d \psi(|X_t-Y_t|)\le & \{C\psi(|X_t-Y_t|) \d t +\theta_t\|\psi'\|_\infty\W_{\psi} (P_t^*\mu, P_t^*\nu) \big\}\d t\\
& + \psi'(|X_t-Y_t|) \<\{\hat\si_t(X_t)-\hat\si_t(Y_t)\}^* u(X_t,Y_t), \d W_t^{(2)}\>.\end{align*}
Noting that $\psi(|X_\tau-Y_\tau|)=0$, we obtain
$$\E\big[1_{\{t>\tau\}}\psi(|X_t-Y_t|)\big] \le \|\psi'\|_\infty\e^{C(t-s)}  \E \int_{t\land\tau}^t\theta_r  \W_{\psi}(P_r^*\mu, P_r^*\nu) \d r,
\ \ t\ge s.$$
Combining this with \eqref{W1'} and \eqref{O1}, we derive
\beg{align*} \W_\psi(P_t^*\mu,P_t^*\nu)&\le \E\psi(|X_t-Y_t|)  = \E\psi(|X_{t\land\tau}- Y_{t\land \tau}|) + \E\big[1_{\{t>\tau\}}\psi(|X_t-Y_t|)\big]\\
 &\le \e^{-\int_s^t\zeta_r\d r}\E \psi(|X_s-Y_s|)  +  \|\psi'\|_\infty \e^{C(t-s)}  \int_s^t  \theta_r\W_{\psi}(P_r^*\mu, P_r^*\nu) \d r\\
&=  \e^{-\int_s^t\zeta_r\d r}\W_{\psi} (P_s^*\mu, P_s^*\nu)   + \|\psi'\|_\infty  \e^{C(t-s)}   \int_s^t \theta_r \W_{\psi}(P_r^*\mu, P_r^*\nu) \d r,
 \ \ t\ge s.\end{align*}
 Therefore,
 \beg{align*} \ff{\d^+}{\d s} \W_{\psi} (P_s^*\mu, P_s^*\nu)&:= \limsup_{t\downarrow s} \ff{\W_{\psi} (P_t^*\mu, P_t^*\nu) -\W_{\psi} (P_s^*\mu, P_s^*\nu) }{t-s} \\
 &\le -(\zeta_s -\theta_s \|\psi'\|_\infty) \W_{\psi} (P_s^*\mu, P_s^*\nu),\ \ s\ge 0.\end{align*}
 This implies \eqref{EXP1'}.
\end{proof}

 As a consequence of Theorem \ref{T7},   we consider the non-dissipative case where $\nn b_t(\cdot,\mu)(x)$ is  positive definite
in a possibly unbounded set but with bounded  $``$one-dimensional puncture mass" in the sense of  \eqref{H23} below.

Let    $\W_1=\W_\psi$ and $\scr P_1(\bar D)= \scr P_\psi$ for $\psi(r)=r$, and define
 \beg{align*}  S_b(x):= \sup\big\{\<\nn_v b_t(\cdot,\mu)(x), v\>:\ t\ge 0, |v|\le 1, \mu\in \scr P_1(\bar D)\big\},\ \ x\in \bar D.\end{align*}

 \beg{enumerate} \item[{\bf(A8)}](3)   There exist constants $\theta_0,\theta_1,\theta_2,\bb \ge 0$     such that
\beq\label{H21}  \ff 1 2 \|\si_t(x)-\si_t(y)\|_{HS}^2\le \theta_0|x-y|^2,\ \ t\ge 0, x,y\in  \bar D;\end{equation}
\beq\label{H22} S_b(x)\le\theta_1, \ \ |b_t(x,\mu)-b_t(x,\nu)|\le \bb \W_1(\mu,\nu),\ \ t\ge 0, x\in\bar D, \mu,\nu\in \scr P_1(\bar D);\end{equation}
\beq\label{H23} \zeta:=  \sup_{x,v\in \bar D, |v|=1}\int_\R 1_{\{S_b(x+sv)>-\theta_2\}}\d s <\infty.\end{equation}
 \end{enumerate}
According to the proof of \cite[Corollary 3.2]{W21a}, the following result follows from  Theorem \ref{T7}.

\beg{cor} \label{NC3.2} Let $D$ be convex. Assume {\bf (A8)}$(1)$ and {\bf (A8)}$(3)$. Let
\beq\label{G2} \beg{split}& \gg(r):= (\theta_1+\theta_2) \big\{(\zeta r^{-1})\land r\big\}-(\theta_2-\theta_0)r,\ \ r\ge 0,\\
&k:= \ff{2\aa}{\int_0^\infty t\, \e^{\ff 1 {2\bb} \int_0^t \gg(u)\d u}\d t }-\ff{\bb (\theta_2-\theta_0)}{2\aa}\int_0^\infty t \e^{\ff 1 {2\aa}\int_0^t \gg(u)\d u}\d t. \end{split}\end{equation}
Then there exists a constant $c>0$ such that
$$\W_1(P_t^*\mu,P_t^*\nu) \le c\e^{-k t}\, \W_1(\mu,\nu),\ \ t\ge 0, \mu,\nu\in \scr P_1(\bar D).$$
If $(b_t,\si_t)$ does not depend on $t$ and $\theta_2>\theta_0$ with
$$ \bb <  \ff{4\aa^2}{(\theta_2-\theta_0)(\int_0^\infty t\, \e^{\ff 1 {2\aa} \int_0^t \gg(u)\d u}\d t)^2 },$$
then  $k>0$ and $P_t^*$ has a unique invariant probability measure $\bar\mu\in \scr P_1(\bar D)$ satisfying
$$\W_1(P_t^*\mu,\bar\mu) \le c\e^{-k t}\, \W_1(\mu,\bar\mu),\ \ t\ge 0, \mu\in \scr P_1(\bar D).$$
\end{cor}

\paragraph{Remark 4.1}  We note that \cite[Theorem 2.1]{W21a} presents an ergodicity result for  the non-dissipative case, which also holds for present setting with convex $D$. We drop the detailed statement to save space.

\paragraph{Acknowledgement.} The  author would like to thank  the referee for helpful comments and corrections.

\end{document}